\documentclass[12pt,a4paper]{amsart}
\theoremstyle{plain}

\usepackage{enumerate, amssymb}

\advance\hoffset-5mm \advance\textwidth40mm
\input{diagrams.tex}
\diagramstyle[scriptlabels,height=8mm,width=8mm]

\def\bdi{\begin{diagram}}
\def\edi{\end{diagram}}

\theoremstyle{plain}

\newtheorem{thm}{Theorem}[section]
\newtheorem{cor}[thm]{Corollary}
\newtheorem{lem}[thm]{Lemma}
\newtheorem{prop}[thm]{Proposition}
\theoremstyle{definition}
\newtheorem{defi}[thm]{Definition}
\newtheorem{defis}[thm]{Definitions}
\newtheorem{conj}[thm]{Conjecture}
\newtheorem{conv}[thm]{Convention}
\newtheorem{nota}[thm]{Notation}
\newtheorem{rem}[thm]{Remark}
\newtheorem{rems}[thm]{Remarks}
\newtheorem{exa}[thm]{Example}
\newtheorem{exas}[thm]{Examples}
\newtheorem{prob}[thm]{Problem}
\newtheorem{probs}[thm]{Problems}
\newtheorem{ques}[thm]{Question}
\newtheorem{sit}[thm]{}

\newcommand{\Ker}{ \operatorname{{\rm Ker}}}

\newcommand{\Hol}{ \operatorname{{\rm Hol}}}

\renewcommand{\epsilon}{\varepsilon}

\def\and{\quad\mbox{and}\quad}

\newcommand{\C}{\ensuremath{\mathbb{C}}}

\newcommand{\Q}{\ensuremath{\mathbb{Q}}}
\newcommand{\Z}{\ensuremath{\mathbb{Z}}}

\newcommand{\proj}{\ensuremath{\mathbb{P}}}
\newcommand{\bk}{{\ensuremath{\rm \bf k}}}

\newcommand{\hA}{{\hat A}}
\newcommand{\hC}{{\hat C}}

\newcommand{\hq}{{\hat q}}
\newcommand{\hx}{{\hat x}}
\newcommand{\hy}{{\hat y}}
\newcommand{\hB}{{\hat B}}
\newcommand{\hD}{{\hat D}}
\newcommand{\hE}{{\hat E}}
\newcommand{\hK}{{\hat K}}
\newcommand{\hL}{{\hat L}}
\newcommand{\hP}{{\hat P}}
\newcommand{\hQ}{{\hat Q}}
\newcommand{\hR}{{\hat R}}
\newcommand{\hS}{{\hat S}}

\newcommand{\hU}{{\hat U}}
\newcommand{\hV}{{\hat V}}
\newcommand{\hW}{{\hat W}}

\newcommand{\hX}{{\hat X}}
\newcommand{\bY}{{\bar Y}}

\newcommand{\tK}{{\tilde K}}
\newcommand{\tQ}{{\tilde Q}}
\newcommand{\tR}{{\tilde R}}
\newcommand{\tS}{{\tilde S}}

\newcommand{\tU}{{\tilde U}}
\newcommand{\tX}{{\tilde X}}
\newcommand{\tY}{{\tilde Y}}
\newcommand{\tZ}{{\tilde Z}}

\newcommand{\bE}{{\bar E}}

\newcommand{\bZ}{{\bar Z}}

\def\bk{{\bf k}}

\newcommand{\tcO}{{\ensuremath{\tilde {\mathcal{O}}}}}

\newcommand{\cB}{{\ensuremath{\mathcal{B}}}}

\newcommand{\cS}{{\ensuremath{\mathcal{S}}}}

\newcommand{\cO}{{\ensuremath{\mathcal{O}}}}

\renewcommand{\rho}{\varrho}

\def\bals#1\eals{\begin{align*}#1\end{align*}}
\def\bal#1\eal{\begin{align}#1\end{align}}

\def\AA{{\mathbb A}}

\renewcommand{\phi}{\varphi}

\newcommand{\bnum}{\begin{enumerate}}
\newcommand{\enum}{\end{enumerate}}
\renewcommand{\emptyset}{\varnothing}

\newcommand{\brem}{\begin{rem}}
\newcommand{\brems}{\begin{rems}}
\newcommand{\erem}{\end{rem}}
\newcommand{\erems}{\end{rems}}
\newcommand{\bprob}{\begin{prob}}
\newcommand{\eprob}{\end{prob}}
\newcommand{\bprobs}{\begin{probs}}
\newcommand{\eprobs}{\end{probs}}
\newcommand{\bques}{\begin{ques}}
\newcommand{\eques}{\end{ques}}
\newcommand{\bexa}{\begin{exa}}
\newcommand{\bexas}{\begin{exas}}
\newcommand{\eexa}{\end{exa}}
\newcommand{\eexas}{\end{exas}}
\newcommand{\bdefi}{\begin{defi}}
\newcommand{\edefi}{\end{defi}}
\newcommand{\bdefis}{\begin{defis}}
\newcommand{\edefis}{\end{defis}}
\newcommand{\bcor}{\begin{cor}}
\newcommand{\ecor}{\end{cor}}
\newcommand{\blem}{\begin{lem}}
\newcommand{\elem}{\end{lem}}
\newcommand{\bconv}{\begin{conv}}
\newcommand{\econv}{\end{conv}}
\newcommand{\bconj}{\begin{conj}}
\newcommand{\econj}{\end{conj}}
\newcommand{\bprop}{\begin{prop}}
\newcommand{\eprop}{\end{prop}}
\newcommand{\bthm}{\begin{thm}}
\newcommand{\ethm}{\end{thm}}
\newcommand{\bnota}{\begin{nota}}
\newcommand{\enota}{\end{nota}}
\newcommand{\bsit}{\begin{sit}}
\newcommand{\esit}{\end{sit}}
\newcommand{\be}{\begin{equation}}
\newcommand{\ee}{\end{equation}}
\newcommand{\bproof}{\begin{proof}}
\newcommand{\eproof}{\end{proof}}
\def\ba{\begin{array}}
\def\ea{\end{array}}

\thanks{ }

\begin{document}

\title[Free $G_a$-actions on $\C^n]{Free $G_a$-actions on $\C^n$}

\author{S.\ Kaliman}
\address{Department of Mathematics,
University of Miami, Coral Gables, FL 33124, USA}
\email{kaliman@math.miami.edu}

%\begin{abstract}  We derive a criterion for a free $G_a$-action on $\C^n$ to be a translation.
%We also prove that every proper $G_a$-action on $\C^4$ preserving a coordinate is a translation.
%\end{abstract}

%\date{\today}
%\maketitle

%\thanks{
%\mbox{\hspace{11pt}}{\it 2000 Mathematics Subject Classification}:
%14R20, 14L30.\\
%{\renewcommand{\thefootnote}{} \footnotetext{ 2010}
%\textit{Mathematics Subject Classification:}
%14R20,\,32M17.\mbox{\hspace{11pt}}\\{\it Key words}: affine\documentstyle[12pt]{report}

%varieties, group actions, one-parameter subgroups, transitivity.}}

%\vfuzz=2pt
%\thanks{}

%
%{\footnotesize \tableofcontents}

\begin{center}{\Large Proper $G_a$-actions on $\C^4$ preserving a coordinate\\[3ex]}\end{center}

\begin{center}{Shulim Kaliman\\[3ex]}\end{center}

\hspace{9cm} {\em Dedicated to my teachers\\

\hspace{9cm}   Vladimir Yakovlevich Lin

\hspace{10.5cm} and 

\hspace{9cm}  Evgeniy Alekseevich Gorin}\\

\vspace{0.5cm}

\begin{center}  ABSTRACT. We prove that the actions mentioned in the title are translations.
We show also that for certain $G_a$-actions on affine fourfolds the categorical quotient of the action is automatically an affine algebraic variety
and describe the geometric structure of such quotients. \\[5ex]
\end{center}

\vspace{1cm}

\section*{Introduction}
An algebraic action $G_a$ of the additive group $\C_+$ of complex numbers on a complex algebraic variety $X$
is free if it has no fixed points. When $X$ is a Euclidean space $\C^n$ with a coordinate system $(x_1, \ldots , x_n)$ the
simplest example of such an action is a translation for which the action of an element $t \in \C_+$ is given by
$(x_1, x_2, \ldots , x_n) \to (x_1+t, x_2, \ldots , x_n)$.  It turns out that for $n\leq 3$ these notions are ``essentially"
the same. More precisely, when $n\leq 3$ every nontrivial free $G_a$-action on $\C^n$ in a suitable polynomial
coordinate system is a translation\footnote{In fact,  for $n\leq 3$ every connected one-dimensional unipotent algebraic subgroup of 
Cremona group of $\C^n$ is conjugate to such a translation \cite[Corollary 5]{Po2}.}  (see \cite{Gu} and \cite{Re} for $n=2$, and \cite{Ka04} for $n=3$).

Starting with $n=4$ the similar statement does not hold and the basic example of Winkelmann \cite{Wi}
gives a triangular\footnote{Recall that a $G_a$-action on $\C^n$ is triangular if in a suitable polynomial coordinate system it is of form
$(x_1, \ldots , x_n) \to (x_1, x_2+ tp_2(x_1), x_3 + tp_3(x_1, x_2), \ldots , x_n+ tp_n(x_1, \ldots x_{n-1}))$ where each $p_i$ is a polynomial.
For $n\geq 3$ not every $G_a$-action on $\C^n$ is triangulable (i.e. triangular in a suitable polynomial coordinate system) \cite{Ba}.}
free $G_a$-action which is not a translation. In his example the geometric quotient of the
action is not Hausdorff while for a translation on $\C^n$ the geometric quotient  % (as well as the algebraic quotient)
is isomorphic to $\C^{n-1}$. Recently Dubouloz, Finston, and Jaradat \cite{DuFiJ} proved that every triangular action on $\C^n$, which is proper
(in particular, it is free and has a Hausdoff geometric quotient), is a translation in a suitable coordinate system.  
Note that every triangular action preserves at least one of coordinates, and one
of the aims of this paper is the following generalization of the Finston-Dubouloz-Jaradat result.

\bthm\label{0.1} Every proper $G_a$-action on $\C^4$ that preserves a coordinate is a translation in a suitable polynomial
coordinate system.\footnote{Neena Gupta informed the author that she and S. M. Bhatwadekar had obtained recently an independent proof of Theorem \ref{0.1}.}
\ethm

Its proof involves investigation of categorical quotients of $G_a$-actions on $\C^4$. Namely, let $X$ be an affine
algebraic variety equipped with a $G_a$-action $\Phi : G_a\times X \to X$.  Denote by $\C [X]^{\Phi}$ the subring of $\Phi$-invariant regular functions
in the ring $\C [X]$ of regular functions on $X$ and by $ {\rm Spec} \, \C [X] ^\Phi$ (resp. $ {\rm Spec} \, \C [X] $) the spectrum of  $\C [X]^{\Phi}$ (resp $\C [X]$). 
Then  the natural embedding $\C [X]^{\Phi}\hookrightarrow \C [X]$ induces a map $ {\rm Spec} \, \C [X]  \to {\rm Spec} \, \C [X]^\Phi $.
The fact that $\C [X]^{\Phi}$
is  finitely generated is equivalent to the fact that $ {\rm Spec} \, \C [X] ^\Phi$ can be viewed as an affine algebraic variety denoted by $X//\Phi$.
It is called the categorical quotient of the action and the map of the spectra yields
the quotient morphism $\rho : X \to X//\Phi$ in the category of affine algebraic varieties.
If $\dim X \leq 3$ then  $\C [X]^{\Phi}$
is  always finitely generated  by a theorem of Zariski \cite{Za}.  
In higher dimensions this fact is not necessarily true by Nagata's counterexample to the fourteenth Hilbert problem.
Furthermore,  extending Nagata's counterexample Daigle and Freudenburg  \cite{DaFr} showed for a $G_a$-action on $\C^n$ the ring of invariant functions  may not be finitely generated starting from dimension $n \geq 5$.
In dimension 4 the same authors
showed that the ring of regular functions invariant with respect to a triangular $G_a$-action on $\C^4$ is automatically finitely generated \cite{DaFre} and later Bhatwadekar 
and Daigle \cite{BhDa} proved that it remains finitely generated
if one considers instead of triangular $G_a$-actions the wider class of $G_a$-actions preserving a coordinate\footnote{The author is grateful to Neena Gupta
for drawing his attention to the paper of  Bhatwadekar 
and Daigle.}. In this paper we establish a stronger fact
contained in the next theorem together with a generalization of Theorem \ref{0.1}.

\bthm\label{0.2} Let $\varphi : X \to B$ be a surjective morphism of a factorial affine algebraic $G_a$-variety $X$ (i.e.  $X$ is equipped 
with some $G_a$-action $\Phi$ ) of dimension 4
into a smooth affine curve $B$.  %Suppose that $\rho : X \to X//\Phi =:Q$ is  the categorical quotient morphism in the category of affine schemes.     
%(i.e. $Q$
%is the spectrum of the subring of $G_a$-invariant functions in the ring of regular functions $\C [X]$ on $X$).
Suppose also that

$\bullet$ the action preserves each fiber of $\varphi$;

$\bullet$  the generic fiber of $\varphi$ is a three-dimensional variety $Y$ (over the field $K$
of rational functions on $B$)
%$\bullet$ for the algebraic closure $\tK$ of $K$ and the variety $\tY$ over $\tK$ obtained from $Y$
for which the ring of invariants of the $G_a$-action (induced by $\Phi$) on $Y$ is
the polynomial ring $K [z,w]$ in two variables; \footnote{This assumption that the ring of invariants of the induced action is isomorphic to $K[z,w]$ can be replaced
by the following: for a general $b \in B$ the categorical quotient of the restriction of $\Phi$ to   the fiber $X_b=\varphi^{-1} (b)$ is isomorphic to $\C^2$ (see,  Theorem \ref{ref.t4} below).}

$\bullet$ for every $b \in B$ the fiber $X_b=\varphi^{-1} (b)$ admits a nonconstant morphism into a curve
if and only if this curve is a polynomial one (i.e. the normalization of the curve is the line $\C$).

%$\bullet$ each fiber $X_b$ is  topologically unibranch at every point $x \in X_b$ (i.e. the analytic germ of $X_b$ at $x$ is irreducible).
Then 

{\rm (1)}  the ring of $\Phi$-invariant functions is finitely generated and, thus, it can be viewed as the ring of regular functions on an affine algebraic variety $Q=X//\Phi$;

{\rm (2)} there is an affine 
modification $\psi : Q \to B\times \C^2$ such that for some nonempty Zariski dense subset $B^* \subset B$
the restriction of $\psi$ over $B^*$ is an isomorphism and every singular fiber of $\psi$ is of form $C \times \C$ 
where $C$ is a polynomial curve.

{\rm (3) } Furthermore, if one requires additionally that 

$\bullet$ $\Phi$ is proper and $X$ is Cohen-Macaulay,

$\bullet$ each fiber $X_b$ is  normal, %locally factorial, 

$\bullet$ and the restriction $\Phi_b$ of $\Phi$ to $X_b$ is a translation (in particular, $X_b$ is naturally isomorphic to a direct product $(X_b//\Phi_b) \times \C$), 

\noindent then the quotient $Q$ is locally trivial $\C^2$-bundle over $B$ (and in particular it is a vector bundle by \cite{BCW}),
$X$ is naturally isomorphic to $Q\times \C$, and $\Phi$ is generated by a translation on the second factor of $X\simeq Q \times \C$.

\ethm

Let us emphasize that the fact that $X$ is a direct product $Q\times \C$ is a rather rare event for regular $G_a$-actions
while in the category of rational actions of a connected linear algebraic 
group $G$ on an algebraic variety $Y$ (over an algebraically closed field of any characteristic) this variety is automatically birationally isomorphic
to the product of $\proj^s$ and the rational quotient of $Y$ with respect to a Borel subgroup of $G$ (see \cite{Ma63} and \cite{Po1}).
It is also worth mentioning that the requirement that the restriction of $\Phi$ to any $X_b$ is a translation in (3) can be can be replaced 
by some topological assumptions. For instance, if each fiber $X_b$ is smooth and factorial with trivial second and third homology groups
then $\Phi|_{X_b}$ is automatically a translation by \cite[Theorem 5.1]{Ka04}. In particular, this is true when $X_b$ is a
smooth contractible threefold since by the Gurjar theorem \cite[Theorem 1]{Gur} (see also \cite{Fu}) such a threefold is factorial.

This leads to the following application of Theorem \ref{0.2}.

\bcor\label{0.3} Let $\varphi : X \to B$ be a surjective morphism of a smooth factorial affine algebraic fourfold $X$
into a smooth curve $B$ such that  $X$ is equipped with a proper $G_a$-action $\Phi$  preserving every fiber $X_b=\varphi^{-1}(b), \, b \in B$
of $\varphi$. Suppose that each $X_b$ is a smooth
contractible threefold such that the restriction of $\Phi$ to $X_b$ has the plane as the quotient \footnote{For $X_b=\C^3$ the quotient
of a nontrivial $G_a$-action is always isomorphic to $\C^2$
due to Miyanishi's theorem \cite{Miy} and this is also true for smooth contractible threefolds
with negative logarithmic Kodaira dimension \cite{KaSa} like the Russell cubic $\{ x+x^2y+z^2+t^3=0 \} \subset \C^4_{x,y,z,t}$.}. 

Then there exits a categorical quotient of the action $Q=X//\Phi$ in the category of affine algebraic varieties. Furthermore, $Q$
is a two-dimensional vector bundle over $B$ and  $X$ is naturally isomorphic to $Q\times \C$ while
$\Phi$ is generated by a translation on the second factor of $X\simeq Q \times \C$.
\ecor

Indeed, it is straightforward that the assumptions of Corollary \ref{0.3} imply all assumptions of Theorem \ref{0.2} with
a possible exception of the condition on the generic fiber of $\varphi$. But this last condition follows from a combination
of the Kambayashi \cite{Kam} and the Kraft-Russell \cite{KrRu} theorems (see Theorem \ref{4.1} and Example \ref{4.3} below for details).

When $X\simeq \C^4$, $B\simeq \C$, and $\varphi$ is a coordinate function the assumptions
of Corollary \ref{0.3} are automatically true and we have Theorem \ref{0.1}. 

%The Dubouloz-Finston-Jaradat  theorem uses in an essential way the important result due to Daigle and Fruedenburg
%\cite{DaFre} who showed that the quotient of a triangular $G_a$-action on $\C^4$ is automatically affine \footnote{In 
%general the algebraic quotient $\C^n//G_a$ is not affine due to Nagata's counterexample to the fourteenth Hilbert problem.
%Furthermore, it may not be affine starting from dimension $n \geq 5$ \cite{DaFr} due to Daigle and Freudenburg,
%while for $n \leq 3$ it is affine by a theorem of Zariski \cite{Za}.}.
%Hence the first statement of Theorem \ref{0.2} strengthens the Daigle-Freudenburg result.

Some of the results mentioned before  (including Theorem \ref{0.1}) are extended in the last section of this paper where 
using the Lefschetz principle we show that similar facts (see, Theorems \ref{7.0.1} and \ref{7.0.2}) remain valid if we consider varieties $X$ and $B$ not over the field of
 complex numbers $\C$ but over any (not necessarily algebraically closed) field of characteristic zero.\\

{\em Acknowledgements.} The author would like to thank the referee for his contribution to the quality of this paper, for his patience, and
ability to see through inaccuracies and mistakes in the original version of this manuscript.

\section{General facts}

\bthm\label{1.1} Let $\varphi : X \to Y$ be a birational morphism of irreducible affine algebraic varieties such that $Y$ is normal
and there is a subvariety $Z$ of $Y$ with codimension at least 2 for which 
the restriction of $\varphi$ to $X\setminus \varphi^{-1}(Z)$ is a surjective morphism onto $Y\setminus Z$.
%$Y \setminus Z \subset \varphi (X)$,
Suppose also that for every point $y \in Y \setminus Z$ the preimage $\varphi^{-1} (y)$ is finite.
%, and
%for a general point $y \in Y\setminus Z$ the preimage $\varphi^{-1} (y)$ is a singleton.
Then $\varphi : X \to Y$ is an isomorphism.
\ethm

\bproof %By assumption $\varphi$ is a birational map. Hence t
The Zariski Main theorem (e.g., see \cite[Th\'eor\`eme 8.12.6]{EGA} or \cite[Chapter III, Corollary 11.4]{Har}) implies that 
the restriction of $\varphi$ yields an isomorphism between $X \setminus \varphi^{-1}(Z)$ and $Y \setminus Z$.
By the Hartogs theorem the composition of the inverse map $\varphi^{-1} : Y \setminus Z \to X \setminus \varphi^{-1} (Z)$
with the inclusion $X \setminus \varphi^{-1}(Z) \hookrightarrow X$
extends to a morphism $Y\to X$ and we are done. 

\eproof

\bcor\label{1.3}   Let $\varphi : X \to Y$ be a morphism of irreducible affine algebraic varieties, 
and $E\subset X$ and $D \subset Y$ be irreducible divisors such that the restriction
of $\varphi$ yields an isomorphism $X \setminus E \to Y \setminus D$.
Suppose also that $Y$ is normal and $\varphi (E)$ is Zariski dense in $D$.
Then $\varphi : X \to Y$ is an isomorphism.

\ecor

\bproof The dimension argument implies that for a general point $y$ in $D$ the preimage
$\varphi^{-1}(y)$ is finite. Hence it is finite outside a proper closed subvariety $Z$ of $D$
and we are done by Theorem \ref{1.1}.

\eproof

The next fact and the modified proof of Theorem \ref{1.4} below were suggested by the referee.

\bprop\label{ref.p1}  Let $\rho : X \to Q$ be a morphism of irreducible affine algebraic varieties such that $Q$ is normal and $R=Q \setminus \rho (X)$ is of codimension
at least 2 in $Q$. Then every regular function $f$ on $X$ which is constant on the general fibers of $\rho$ descends to a unique regular function $g$ on $Q$.

\eprop

\bproof 

Let us show first that $f$ is a lift of a rational function $g$ on $Q$.
That is, one needs to establish the regularity
of $g$ on a Zariski dense open subset $Q_0$ of $Q \setminus R$.  We can suppose that the restriction of $\rho$ to $X_0=\rho^{-1}(Q_0)$ is flat
by \cite[Th\'eor\`eme 6.9.1]{EGA} and, therefore, being surjective it is faithfully flat \cite[Ch. 3, Exercise 16]{AM}. 
%[Theorem ?]{Mat}. % \cite[Th\'eor\`eme 6.7.8]{EGA01}. % Note also that $\rho|_{X_0} : X_0 \to Q_0$ can be viewed as an fpqc cover (e.g., 
%see \cite{Vis}) for terminology).

%Consider $P =X_0 \times \C$,  the graph $\Gamma \subset P$ of $f|_{X_0}$, and $T=\Gamma \times_{Q_0} \Gamma \subset P \times_Q P$.
%Let $\pr_i : T \to \Gamma \subset X_0 \times \C, \, i=1,2$ be the natural projections.
%Since $f$ is constant on every fiber $\rho^{-1} (q), \, q \in Q_0$ we see that $\pr_1(T)=\pr_2 (T)$. By \cite[Corollary 4.34]{FGIKNV} % (or \cite[Theorem 2]{FGA}) 
%this implies that $\Gamma$ is a lift of a subvariety $\Lambda \subset Q_0 \times \C$ via the morphism induced by $\rho$. This subvariety $\Lambda$ yields the graph of the regular function $g|_{Q_0}$ 
%Hence $X_0$ can be viewed as an fpqc covering of $Q_0$ and $f|_{X_0}$ can be viewed as a section of an  quasi-coherent $\cO_{Q_0}$-module which is a descent data
%with respect to $X_0\to Q_0$.  Since the functor from the category of  quasi-coherent $\cO_{Q_0}$-modules on $Q_0$
%to the category of descent data with respect to such a covering is fully faithful (e.g., see \cite[Chapter 34, Lemma 5.1]{Columbia}  or  \cite[Theorem 4.23]{FGIKNV})  $f|_{X_0}$ must be a lift
%of a regular function $g|_{Q_0}$ which implies that $g$ is rational on $Q$.
Consider $Y=X_0\times_{Q_0}X_0$. Then the two natural projections $Y \to X_0$ generate two homomorphisms $e_1:\C[X_0] \to \C[Y]$ and $e_2: \C[X_0] \to \C[Y]$.
Since $f$ is constant on general fibers of $\rho$ we see that $f|_{X_0} \in \Ker (e_1-e_2)$. In combination with the faithful flatness of the natural morphism $\C [Q_0] \to \C[X_0]$ this implies that $f$ must be a lift
of a regular function $g|_{Q_0}$  \cite[Lemma 2.61]{FGIKNV}. Hence  $g$ is a rational function on $Q$.

Assume that there exists a divisor $T$ in $Q$ such that  a general point $t_0 \in T$ does not belong to the indeterminacy set of $g$ and $g(t_0)=\infty$. 
We can suppose that also $t_0 \notin R$, i.e. $\rho^{-1}(t_0) \ne \emptyset$. Hence for a germ of a curve $C$ in $X$ through $x_0 \in \rho^{-1}(t_0)$ one has
$f(c) \to f(x_0)$ as $c \in C$ approaches $x_0$. This implies that $g(\rho (c))\to f(x_0) \ne \infty$ as $\rho (c)$ approaches $t_0$. A contradiction.
That is, $g$ is regular on $Q$ outside a subvariety of  codimension at least 2. Since $Q$ is normal the Hartogs'  theorem implies that $g$ is regular on $Q$
and we are done. 

\eproof

\bthm\label{1.4}  Let $\rho : X \to Q$ be a morphism of irreducible affine algebraic varieties and $\Phi : G_a \times X\to X$ be a non-trivial $G_a$-action
on $X$ which preserves each fiber of $\rho$.
Let $Q$ be normal and $P$ be a subvariety of $Q$ with codimension at least 2 such that
$Q\setminus P \subset \rho (X)$. Suppose also that
for every point $q \in Q \setminus P$ the preimage $\rho^{-1} (q)$ is a curve and
for a general point $q \in Q\setminus P$ the preimage $\rho^{-1} (q)$ is an irreducible curve.
Then $\rho : X \to Q$ is  the categorical quotient morphism in the category of affine algebraic varieties.  In particular, the subring of $G_a$-invariants
in the ring of regular functions on $X$ is finitely generated.
\ethm

\bproof Note that by the assumption every general fiber of $\rho$ is nothing  but an orbit of $\Phi$. Thus every $\Phi$-invariant function $f$ is constant on the general fibers of $\rho$.
By Propostion \ref{ref.p1} $f$ is a lift of a regular function on $Q$. Hence the ring of regular functions on $Q$
coincides with the ring  $\C [X]^{\Phi}$  of invariants and we are done.

%Suppose that $f$ is a $\Phi$-invariant regular function on $X$. Consider the map
%$\rho_1 = (\rho , f): X \to Q \times \C$. Let $Q_1$ be the closure of its image and $\tau :Q_1 \to Q$
%be the natural projection. Since the image of an irreducible variety is irreducible, so is $Q_1$. 
%Furthermore, since the general fibers of $\rho$ are irreducible curves, the morphism $\tau$ is birational.
%This implies that for a subvariety $P' \supset P$ of $Q$ which codimension is at least 2
%the restriction of $\tau$ to $\tau^{-1}(Q \setminus P')$ is quasi-finite (indeed, the minimal $P'$
%over which the restriction of $\tau$ is not quasi-finite
%cannot be a divisor since otherwise $Q_1$ is not irreducible).
%Theorem \ref{1.1} implies that  
%$\tau : Q_1 \to Q$ is an isomorphism. Hence $f$ can be viewed as a regular function on
%$Q$. Since $f$ is an arbitrary $\Phi$-invariant function the ring of regular functions on $Q$
%coincides with the whole ring  $\C [X]^{\Phi}$  of invariants and we are done.

\eproof

Now the dimension argument as in the proof of Corollary \ref{1.3} implies the following.

\bcor\label{1.5}   Let $\rho : X \to Q$ be a morphism of irreducible affine algebraic varieties such that $Q$ is normal
and let  $\Phi : G_a \times X \to X$  be a non-trivial $G_a$-action on $X$ preserving every fiber of  $\rho$. 
Let $E\subset X$ and $D \subset Q$ be irreducible divisors such that  $X\setminus E$ and $Q \setminus D$ are affine algebraic varieties for which
$\rho (X \setminus E)= Q\setminus D$ and $\rho |_{X\setminus E)} : X \setminus E \to Q \setminus D$
is the categorical quotient morphism of the action $\Phi |_{X \setminus E}$.
Suppose also that $\rho (E)$ is Zariski dense in $D$.

Then $\rho : X \to Q$ is the categorical quotient morphism of $\Phi$. 

\ecor

\bprop\label{fv.p1} {\rm (cf. \cite[Lemma 2.1]{Ka04})} Let $\rho : X \to Q$ be a dominant morphism of normal affine algebraic  varieties.
Suppose that the general fibers of $\rho$ are irreducible and there are no non-constant invertible regular function on such fibers. 
Suppose also that $Q \setminus \rho (X)$ is of codimension at least 2 in $Q$.
Then

%{\rm (1)} for every $f \in \C [X]$, that is constant on each general fiber of $\rho$, one has $f =g\circ \rho$ where $g\in \C [X]$;

{\rm (1)} for every principal irreducible divisor $D$ in $X$, which does not meet general fibers of $\rho$, the closure of $\rho (D)$ is the support of a principal irreducible divisor in $Q$;

{\rm (2)} if $X$ is a factorial variety so is $Q$; 

{\rm (3)}  if $X$ is factorial the preimage of any irreducible reduced  divisor $T$ in $Q$ is an irreducible  reduced divisor in $X$.

\eprop

\bproof  
%Consider  $f$ from (1). 
%By continuity $f$ is constant on every 
%fiber $\rho^{-1}(q)$ whose dimension is $\dim X -\dim Q$.
%Let $P=\{ q \in Q | \dim \rho^{-1}(q) > \dim X - \dim Q \}$. 
%That is, $f$ induces a continuous function $g$ on %an open Zariski dense subset of $\rho (X)\setminus P$ which is locally bounded on  
%$\rho (X)\setminus P$ which is regular because of normality of $Q$. 
%Note that $P$ is of codimension at least 2 in $Q$ since otherwise $X$ is not irreducible.
%Hence $g$ is regular on $Q$ by the Hartogs'  theorem. Thus we have (1).

%{\em Claim.}  For every irreducible divisor $D$ in $X$, which does not meet general fibers of $\rho$, the closure of $\rho (D)$ is the support of an irreducible divisor in $Q$.

In (1) $D$ is the zero locus of a regular function $f$.
 Since $f$ does not vanish on general fibers and they are irreducible we see that $f$ is constant on each general fiber.
By Proposition \ref{ref.p1}  $f =g\circ \rho$ where $g \in \C [Q]$. Let us show that the zero locus of $g$ coincides with the closure of $\rho (D)$. 
Assume to the contrary that there is a divisor $F \subset g^{-1}(0) \setminus \rho (D)$ in $Q$.
Choose a rational function $h$ on $Q$ whose poles are contained in the closure $T$ of  $F$ and a regular function $e$ on $Q$
that vanishes on $T \cap \rho (D)=T \cap \rho (X)$ but not at general points of $T$. Then for sufficiently large $k$ the function $(e^kh)\circ \rho$ is regular on $X$. Since
it is constant on every general fiber of $X$ we conclude by  Proposition \ref{ref.p1} that  $e^kh$ is 
regular on $Q$ contrary to the fact that this function has poles on $T$. Thus we have (1).

Let $T$ be an irreducible reduced Weil divisor in $Q$ and let $D$ be an irreducible component of $\rho^{-1}(T)$ that is a divisor in $X$
(such a component exists because  $Q \setminus \rho (X)$ is of codimension at least 2 in $Q$).
Under the assumption of (2) $D=f^*(0)$ for a regular
function $f$ on $X$ since $X$ is factorial.  By (1) $T=\overline{\rho (D)}$  and it coincides with the zeros of $g \in \C [X]$ where $f=g \circ \rho$.
Furthermore, if $e$ is another regular function on $Q$ that vanishes on 
$T$ then it is divisible by $g$ since $e \circ \rho$ is divisible by $f$, and by  Proposition \ref{ref.p1} $e\circ \rho/g \circ \rho$ is the lift of a regular function on $Q$. That is, 
$T=g^*(0)$. This implies that every irreducible element in $\C [Q]$ has a zero locus which is reduced irreducible and principal.  In particular, this element is prime and,
therefore, we have (2).

Assume that $D_1$ and $D_2$ are reduced irreducible components of $\rho^{-1}(T)$.  Since $X$ is factorial $D_i=f_i^*(0)$.  By the argument above $f_i=g_i\circ \rho$
for regular functions $g_i$ on $Q$.  Furthermore $g_i$ vanishes on $T$ only and thus $g_1/g_2$ is an invertible regular function. This implies that
$D_1=D_2$ and we have (3).

\eproof 

\bcor\label{fv.c1} Let $X$ be a normal affine algebraic variety and $\rho : X \to Q$ be a quotient morphism of a nontrivial $G_a$-action on $X$ (in the category of affine 
algebraic varieties).
 Then $Q \setminus \rho (X)$ has codimension at least 2 in $Q$. Furthermore, if $X$ is factorial so is $Q$.
%Then $Q$ is a factorial variety.
\ecor

\bproof Note that $Q$ is normal since $X$ is. % Furthermore, $Q\setminus \rho (X)$ must be of codimension at least 2 in $Q$.
If $Q\setminus \rho (X)$ contains a divisor $F$ of $Q$ then as in Statement (1)  of Proposition \ref{fv.p1} we can construct a rational function $g$
on $Q$ with poles on $F$ whose lift to $X$ is regular. By construction this lift is $G_a$-invariant, i.e. $g$ must be regular. 
A contradiction.  This implies the first statement while the second one follows  from Proposition \ref{fv.p1}.

\eproof

\section{Basic definitions and properties of affine modifications}

\bdefi\label{ru.d2.1} Recall that an affine modification is any birational morphism $\sigma : \hX \to X$ of affine algebraic varieties \cite{KaZa}. In particular, there exists a divisor
$D \subset X$ such that for $\hD  =\sigma^{-1}(D)$ the restriction of $\sigma$ yields an isomorphism $\hX \setminus \hD \to X\setminus D$.
There is some freedom in the choice of $D$ and we can always suppose that $D$ is a principal effective divisor given  by zeros of a regular function
$f \in A:= \C [X]$. This is the case which we consider in the present paper.  When such $f$ is fixed we call $D$ the divisor of modification, $\hD$
is called the exceptional divisor of modification, and the closure $Z$
of $\sigma (\hD)$ in $X$ is called the center of modification. The advantage of a principal divisor $D$ is that the algebra $\hA =\C [\hX ]$ 
can be viewed as a subalgebra $A [\frac{I}{f}]$ in the field ${\rm Frac} (A)$ of fractions of $A$ where $I$ is an ideal in $A$ (see \cite{KaZa}). 

It is worth mentioning that  $I$ is not determined uniquely, i.e. one can find another ideal $J\subset A$ for which $\hA =A [\frac{J}{f}]$. We suppose further
that $I$ is the largest among such ideals $J$ and we  call $I$ the ideal of the modification. (Treating $A$ as a subalgebra of $\hA$ one can see that
$I$ is the intersection of $A$ and the principal ideal generated by $f$ in $\hA$.)
\edefi

\bnota\label{mod1a} The symbols $X, Z, D, \hX, \hD, I, f$ in this section have the same meaning as in Definition \ref{ru.d2.1}.
%Furthermore throughout the paper for every algebraic variety $X$ by $\C [X]$ we denote its ring of regular functions
%while in the case of a complex space $X$ the ring of holomorphic functions on $X$ will be denoted by $\Hol (X)$.
\enota

\bexa\label{mod2}  Let $I$ be generated by regular functions $f, g_1, \ldots , g_n\in \C [X]$.
Consider  the closed subspace $Y$ of $X \times
\C^n_{v_1, \ldots ,v_n}$ given by the system of equations $fv_j=g_j$ ($j=1, \ldots , n$) and the proper transform
$\hX$ of $X$ under the natural projection $Y \to X$  (i.e. $\hX$ is the only irreducible component of $Y$ whose image in $X$ under the projection is dense in $X$). 
Then the restriction $\sigma : \hX \to X$ of the natural projection
is our affine modification.  Note that $D$ (resp. $\hD$) coincides with the zero locus 
of $f$ (resp. $f \circ \sigma$) and the center $Z$ is given by the equations $f=g_1= \ldots =g_n=0$. 

\eexa

\brem\label{mod3} (1)
The geometrical construction behind the modification in Example \ref{mod2} is the following.  Consider blowing-up $\tau : {\tilde X} \to X$
of $X$ with respect to the ideal sheaf generated by  $f$ and
$g_{1}, \ldots ,g_n$. Delete from ${\tilde X}$ divisors on which
the zero multiplicity of $f \circ \tau$ is more than the zero
multiplicity of at least one of the functions $g_j \circ \tau$. The
resulting variety is $\hX$.

(2) Note that the replacement in  Example \ref{mod2} of functions $ f, g_1, \ldots , g_n$ by functions $ hf, hg_1, \ldots , hg_n$ respectively (where $h \in \C [X]$ is nonzero)
does not change the modification $\sigma : \hX \to X$. In order to avoid this ambiguity we have to fix $f$. We are not going to specify such $f$'s for affine modifications
considered below since the choice of $f$ will be clear from the context in each particular case.

%(2) If we deal with an affine modification in Example \ref{mod2} then the ring $ \C [\hX ]$ of regular function on $\hX$
%is generated by $g_1/f, \ldots , g_n/f$ over $\C [X]$ (see \cite{KaZa}). In particular, $\hX$ is also an affine algebraic
%variety. 

\erem

\bdefi\label{mod4} 
(1)  Let the center of modification $Z$ in Example \ref{mod2} be a set-theoretical complete intersection in $X$
given by the zeros $f=g_1= \ldots =g_n=0$ (i.e. $Z$ coincides with the set of common zeros of these functions and the codimension of $Z$ in $X$ is $n+1$). Then we call such $\sigma : \hX \to X$ a Davis modification.
Its main property is that $\hD$ is naturally isomorphic to $Z \times \C^n$ and that the support of $\hX$ coincides with the support of $Y$ 
\footnote{Indeed, the preimage of $Z$
in $Y$ is naturally isomorphic to $Z\times \C^n$ and therefore has dimension $\dim X - 1$. On the other
hand counting the number of equations defining $Y$ in $X \times \C^n$ we see that every irreducible component of $Y$
has dimension at least $\dim X$. This implies that $Y$ is irreducible (since all irreducible components of $Y$ but one
are contained in the preimage of $Z$) and thus $\hD \simeq Z \times \C^n$.}.

(2) Let $Z$ be a strict complete intersection in $X$ given by $f=g_1= \ldots =g_n=0$; that is, $Z$ is not only set-theoretical 
complete intersection but also the defining ideal of the (reduced) subvariety $Z$ in $X$ coincides with the ideal $I$ generated by $f, g_1, \ldots  , g_n$. Then we call $\sigma$ a simple 
modification. In this case $\hX$ coincides with $Y$ as a scheme.  Note also that
the zero multiplicity of $f \circ \sigma $ is 1 at general points of $\hD$.

\edefi

Here are some useful properties of simple modifications from
\cite{Ka02} which we shall need later.

\bprop\label{normality} Let $\sigma : \hX \to X$ be a simple
modification. Then 

{\rm (1)} $\hX$ is smooth over points from $Z_{\rm reg} \cap D_{\rm reg} \cap X_{\rm reg}$;

{\rm (2)} $\hX$ is Cohen-Macauley provided $X$ is
Cohen-Macauley;

{\rm (3)} furthermore, if in (2)  $X$ is normal and none of irreducible components of the center
$Z$ of $\sigma$ is contained in the singularities of $X$ or $D$ then  $\hX$ is normal. \eprop

\section{Pseudo-affine modifications}

From a geometrical point of view it is sometimes convenient for us to consider a neighborhood $U \subset X$ of some point $z \in Z$ in the standard
topology (we call such $U$ a Euclidean neighborhood) and the restriction $\sigma|_\hU : \hU \to U$ where $\hU$ is any connected component of $\sigma^{-1}(U)$. Since $U$ and $\hU$
are not algebraic varieties but only complex spaces let us consider the analogue of affine modifications in the analytic setting \footnote{One can adhere to
the algebraic setting by viewing $U$ as an \'etale neighborhood of $z$ in $X$.}.

\bdefi\label{mod1} Let $X$ be an irreducible complex Stein space and 
$\psi : X -\to W$ be  a meromorphic map into a projective
algebraic variety $W$ with a fixed ample divisor $H$ such that $\psi$ is holomorphic over $W \setminus H$. 
A minimal resolution   $\pi : \tX \to X$  of indeterminacy points for $\psi$ leads to  a holomorphic map $\tilde \psi : {\tilde X} \to W$.
Removing from ${\tilde X}$ the preimage of $H$ we obtain $\hX$
which, together with the natural projection $\sigma : \hX \to X$,
will be called a pseudo-affine modification.  Consider the Weil divisor $\tilde \psi^*(H)$ and its push-forward (as a cycle ) $D$ by $\pi : \tX \to X$. %The preimage $D=\psi^{-1}(H)\subset X$
Then $D$ is the divisor of
the modification, $\hD = \sigma^{-1} (D)\subset \hX$ is its exceptional divisor, and the closure $Z$ of $\sigma (\hD)$
is its center.  The restriction of $\sigma$ induces, of course, a biholomorphism between $\hX \setminus \hD$ and $X \setminus D$.
Note also that similarly to the algebraic setting $\sigma |_{\hD} :
\hD \to \sigma (\hD )$ is, by construction,  the restriction of the proper morphism $\pi |_{T} : %\tilde \psi^{-1}(H) 
T\to Z$  where $T$ is the closure of $\hD$ in $\tX$. The latter observation is important for the next remark.
 %(this observation will allow us later to use the Stein factorization theorem for $\sigma  |_{\hD}$).
 \edefi

\brem\label{ru.r2.1a}  Though we mostly omit explicit formulations it should be emphasized that
practically all the facts valid for affine modifications have similar analytic analogues for pseudo-affine modifications (e.g. see Proposition \ref{an.normality} below).
\erem

\bexa\label{ru.mod2}  Let us switch to the analytic setting in Example  \ref{mod2}  by assuming that $X$ is a Stein variety, $f, g_1, \ldots ,$ and $g_n$ are holomorphic functions on $X$.
 As before $Y$ is given in $X \times
\C^n_{v_1, \ldots ,v_n}$ by equations $fv_j=g_j$ ($j=1, \ldots , n$) and $\hX$ is the proper transform
of $X$ under the natural projection $Y \to X$.
Then the restriction $\sigma : \hX \to X$ of the natural projection
is  a pseudo-affine modification with $D$ (resp. $\hD$) being the zero locus 
of $f$ (resp. $f \circ \sigma$) and the center $Z$ given by the equations $f=g_1= \ldots =g_n=0$. 

\eexa

\bdefi\label{ru.d2a.1} If in Example \ref{ru.mod2} $Z$ is a set-theoretical complete intersection in $X$ given  by the zeros $f=g_1= \ldots =g_n=0$
then we call $\sigma$ a Davis pseudo-affine modification.
If furthermore $Z$ is a strict complete intersection given by these function then $\sigma$ is a simple pseudo-affine modification.
\edefi

\brem\label{ru.r2.3} 
(1)  For any pair $U, \hU$ as in the beginning of this section the restriction $\hU \to U$ of a simple (resp. Davis) affine modification $\sigma : \hX \to X$  is
automatically a simple (resp. Davis) pseudo-affine modification.

(2) Note that when $X$, $D$, and $Z$ are smooth for a simple pseudo-affine modification then for every point $z \in Z$ the collection
$\{ f, g_1,  \ldots , g_n \}$ can be extended to a local coordinate system in a Euclidean neighborhood $U$ of $z$ in $X$.
In particular, if $n=1$ we can treat $U$ as a germ of $\C^n$ at the origin with coordinates $(u, v, w_1 \ldots , w_{n-2})$
such that $D$ is given in $U$ by $u=0$, $Z$ by $u=v=0$, and the preimage of  $U$ in $\hX$ is viewed as a subvariety of $U\times \C_w$ given
by the equation $uw=v$.

(3) Similarly to the algebraic setting the exceptional divisor of a Davis pseudo-affine modification is a the product of its center and a Euclidean space.
\erem

The following analytic version of Proposition \ref{normality} for simple pseudo-affine modifications remains valid with a verbatim proof.

\bprop\label{an.normality}
Let $\sigma : \hX \to X$ be a simple pseudo-affine 
modification of irreducible Stein spaces. Then 

{\rm (1)} $\hX$ is smooth over smooth points from $Z$ that are not contained in the singularities of  $D$ or $X$;

{\rm (2)} $\hX$ is Cohen-Macauley provided $X$ is
Cohen-Macauley;

{\rm (3)} furthermore, if in (2)  $X$ is normal and none of irreducible components of the center
$Z$ of $\sigma$ is contained in the singularities of $X$ or $D$ then  $\hX$ is normal.

\eprop

\blem\label{coordinates} Let $\psi : Y \to X$ be a holomorphic map of irreducible  normal Stein spaces
such that for some principal effective reduced divisor $D \subset X$ and every $x \in X \setminus D$ the preimage $\psi^{-1}(x)$ is a curve.
Suppose that  the divisor $E=\psi^*(D)$ is reduced and irreducible, $\psi (E)$ is not contained in the singularities of $X$ or $D$,
and for a general point $y\in E$ 
the variety $\psi^{-1} (\psi (y))$ is a surface. Then for such a point $y \in Y$ (resp. for $z=\psi (y) \in X$) there exists a local analytic coordinate
system  $(u',v',v'', \bar w' )$ on $Y$ at $y$  (reps. $(u,v, \bar w )$  on $X$ at $z$) where $\bar w'=(w_1', \ldots , w_{n-2}')$ (resp. $\bar w=(w_1, \ldots , w_{n-2})$)
for which the local coordinate form of $\psi$ is given by
$$(u',v', v'',\bar w') \to  (u,v, \bar w) =(u', (u')^lq(u',v',v'', \bar w'), \bar w')$$
where $l\geq 1$ is the minimal zero multiplicity of the $(n\times n)$-minors in the Jacobi  matrix of $\psi$ at $y$ and the function $q(0,v', v'', \bar w')$ depends on $v'$ or $v''$.
\elem

\bproof  Since $y$ is a general point of $E$ we see that it is a smooth point of $E$ and $X$ and by the assumption $z$ is a smooth point of  $Z =\overline{\psi (E)}, D$, and $X$.
Thus locally $D$ is given by $u=0$ and $Z$ by $u=v=0$ where the holomorphic functions $u$ and $v$ can be included in
a local coordinate system $(u,v, \bar w )$  on $X$ at $z$. If $u'=u \circ \psi$ then by the assumption $E$ is given locally near $y$ by $u'=0$. 
Since $\dim E -\dim Z=2$, by   \cite[ Appendix, Theorem 2]{Chirka}
there exists a local coordinate system $(v',v'', \bar w')$ on $E$ such that the coordinate form of $\psi|_E : E \to D$ is given
by $v=0$ and $\bar w=\bar w'$. Extending functions $v', v''$, and $\bar w'$ holomorphically to a neighborhood of $y$ in $Y$ we get
a local coordinate system $(u',v',v'',\bar w')$ in this neighborhood. Furthermore, the extension $\bar w'$ can be chosen as  $\bar w'=\bar w\circ \psi$.
Hence $\psi$ is given locally by $u=u', \bar w=\bar w'$,  and $v=(u')^lq(u',v',v'', \bar w')$ where $q(0,v', v'', \bar w')$ is not identically zero.
If  $q(0,v', v'', \bar w')$  is independent of $v'$ or $v''$ then
replacing $v$ by $v-u^lq(0, \bar w)$ we increase $l$ (which is at least 1 since otherwise $\psi (E)$ is not contained in $Z$).  On the other hand
$l$ cannot be larger than the minimal zero multiplicity of the  $(n\times n)$-minors in the Jacobi matrix of $\psi$ at $y$. Hence we can suppose
that $q(0,v', v'', \bar w')$ depends on $v'$ or $v''$ in which case $l$ is such a multiplicity.
\eproof

\bprop\label{2.12} 
%Let $\psi : Y \to X$ be a holomorphic map of irreducible  normal Stein analytic sets
%such that for some principal effective reduced divisor $D \subset X$ and every $x \in X \setminus D$ the preimage $\psi^{-1}(x)$ is a curve.
%Suppose that divisor $E=\psi^*(D)$ is reduced and irreducible, $\psi (E)$ is not contained in the singularities of $X$ or $D$,
%and for a general point $y\in E$ 
%the variety $\psi^{-1} (\psi (y))$ is a surface.
Let the assumptions of Lemma \ref{coordinates} hold.
Suppose also that 
\bdi \to X_m & \rTo^{\sigma_m}& X_{m-1}&\rTo^{\sigma_{m-1}}& \ldots &\rTo^{\sigma_2} & X_1 &\rTo^{\sigma_1}& X_0:=X
\edi is a sequence of simple pseudo-affine modifications such that

{\rm (i)} there exists a holomorphic map $\psi_m : Y \to X_m$ for which
$\psi =\sigma_1 \circ \ldots \circ \sigma_m \circ \psi_m$;

{\rm (ii)} for $\psi_i=\sigma_{i+1} \circ \ldots \circ \sigma_m \circ \psi_m$  the closure of  $\psi_i(E)$ coincides with $Z_i\subset D_i$ where
$Z_i$ and $D_i$ are the center and the divisor of $\sigma_{i+1}$ respectively;

{\rm (iii)}  general points of $Z_i$ are contained in the smooth parts of $X_i$ and $D_i$
\footnote{Actually, (iii) follows automatically from the Proposition \ref{an.normality}.   }.

Then such a sequence cannot be extended to the left indefinitely.
% and furthermore for the largest possible $n$ the preimage
%$\psi_n^{-1}(x_n)$ is a curve for every general $x_n \in \psi_n (E)$.
\eprop

\bproof Let $y$ be a general point in $E$ such that $\psi_i(y)=z_i \in Z_i$ which implies that near these points  $Z_i, D_i, E, Y, X$
are smooth. By Lemma \ref{coordinates} we can consider  local coordinate systems $(u',v',v'', \bar w' )$ on $Y$ at $y$ and $(u_i,v_i, \bar w_i )$ on $X_i$ at $z_i$
such that $\psi_i$ is given locally by $u_i=u', \bar w_i=\bar w'$,  and $v_i=(u')^lq(u',v',v'', \bar w')$.
By Remark \ref{ru.r2.3} (2) for a local coordinate system $(u_{i+1}, v_{i+1}, \bar w_{i+1})$ on $X_{i+1}$ 
the coordinate form of $\sigma_{i+1}$ is $(u_{i}, v_{i}, \bar w_{i})=(u_{i+1}, u_{i+1}v_{i+1}, \bar w_{i+1})$.
Hence a local form of $\psi_{i+1}$ is $u_{i+1}=u', \bar w_{i+1}=\bar w'$, and
$v_{i+1}=(u')^{l_i-1}q(u',v',v'',\bar w')$. That is, in this construction $l_{i+1}=l_i-1$. Since such powers cannot be negative
we get the desired conclusion.

\eproof

\brem\label{Jacobi} In fact, we showed that $m$ in Proposition \ref{2.12} cannot exceed
the minimal zero multiplicity $l$ of the $(n\times n)$-minors in the Jacobi matrix of $\psi$ from Lemma \ref{coordinates}.

\erem
%Assume that $Y=\hX \times \C$ and $\psi$ in Proposition \ref{2.12} is the composition of the 
%natural projection and a pseudo-affine modification $\sigma : \hX \to X$. The we get
Similarly, we have the following.

\bprop\label{2.13} Let $\sigma : \hX \to X$ be a pseudo-affine modification 
with center $Z$, divisor $D$, and exceptional divisor $\hD$. Suppose that none of
components of $Z$ is contained in $X_{\rm sing} \cup D_{\rm sing}$, and the image of every
component of $\hD$ is of codimension 1 in $D$. Let 
\bdi \to X_n & \rTo^{\sigma_n}& X_{n-1}&\rTo^{\sigma_{n-1}}& \ldots &\rTo^{\sigma_2} & X_1 &\rTo^{\sigma_1}& X_0:=X
\edi be a sequence of simple modifications with similar conditions on centers and such that $\sigma$ factors through  the composition of these simple modifications.
Then such a sequence cannot be extended to the left indefinitely without violating the fact that $\sigma$ factors through the composition.
\eprop

\bproof For local analytic coordinates systems at general points $y \in \hD$ and $z=\sigma (y) \in Z\subset X$ consider
the zero multiplicity  $k$ of the  Jacobian of $\sigma$ at  $y$ (clearly this multiplicity is independent of
the choice of local coordinates system and the choice of a general point $y$) and let $k_i$ be the similar multiplicity for the modification $\sigma_i \circ \cdots \circ \sigma_1 : X_i \to X$.
Since $\sigma_{i+1}: X_{i+1} \to X_i$ contracts the exceptional divisor in $X_{i+1}$ we see that $k_{i+1}>k_i$. On the other hand if $\sigma$ factors through
$\sigma_i \circ \cdots \circ \sigma_1$ one must have $k_i \leq k$. This yields the desired conclusion.
\eproof

\bdefi\label{ru.d2a.3} (1) Given two (pseudo-) affine modifications $\sigma : \hX \to X$ and $\delta : \tX \to X$ with the same center $Z$ and divisor $D\subset X$
we say that $\sigma$ dominates $\delta$ if it factors through $\delta$. For instance, consider the normalization $\nu : \tX' \to X$ and
let $\delta' = \delta \circ \nu : \tX' \to X$. Then $\delta$ is dominated by $\delta'$.

(2) Let $f,g$, and $h$ be holomorphic functions on a Stein manifold $X$ and 
$\sigma : \hX \to X$ be a simple pseudo-affine modification with a smooth divisor $D=f^*(0)$ and a smooth center $Z$ given by a strict complete intersection
$f=h=0$.
Suppose that $Z$ is also a set-theoretical complete intersection given by $f^k=g=0$ and $\tX \subset X \times \C_w$ is given by $f^kw=g$.
Then we call the Davis modification $\delta : \tX \to X$ (induced by the natural projection) homogeneous  (of degree $k$) if $fI^{k-1}$ and $g$ generate the ideal $I^k$
where $I$ is  the defining ideal  of $Z$ in the ring of holomorphic functions on $X$.
\edefi

\blem\label{ref.l3} Let $X$ be a germ of  $\C^n$ at  the origin  with coordinates $(u, v, w_1 \ldots , w_{n-2})$ and let $\sigma : \hX \to X$ and $I$
be as in Definition \ref{ru.d2a.3} (2) with $f=u$ and $h=v$.  
Suppose also that $g$ and $\delta$ are as in Definition \ref{ru.d2a.3} (2)  without $\delta$ being a priori homogeneous.
% that $g$ is a holomorphic function on $X$ such that the set $Z=\{ u=v=0\}$ coincides also with $\{ u=g=0 \}$.
%Let  $\tX \subset X\times \C_w$ be given by $u^kw=g$ and $\delta : \tX \to X$ be the Davis modification induced by the natural projection.
Then $\delta$ is homogeneous if and only iff the function $g$ is of form
$g=ev^k +a$ where $a\in fI^{k-1}$ %, $b\in I^{k+1}$,
and $e$ is an invertible holomorphic function.

\elem

\bproof  %The ideal $I$ of $\sigma$ is generated by $u$ and $v$ where $u$ plays the role of $f$ in  Definition \ref{ru.d2a.3} (2). 
Since $f=u$ and $I$ is generated by $u$ and $v$ one can see that  $fI^{k-1}$  is generated by functions $u^k, u^{k-1}v, \ldots , uv^{k-1}$.
That is, the ideal $I^k$ is generated by $fI^{k-1}$ and $v^k$ and also by $fI^{k-1}$  and $g$. This is possible iff $g=ev^k +a$ where $a\in fI^{k-1}$
and $e$ is an invertible holomorphic function.
%Since $g\in I^k$ it is of the form $g=ev^k +a$ where $a\in fI^{k-1}$ and $e$ is a holomorphic function. Since $v^k$ is contained in the ideal generated by $g$ and $fI^{k-1}$
%this presentation can be chosen so that $e$ is invertible. 
\eproof

\bprop\label{ru.p2.4} Let $\sigma : \hX \to X$ and $\delta : \tX \to X$ be as in Definition \ref{ru.d2a.3} (2). Then $\hX$ is a normalization of $\tX$ and in particular
$\sigma$ dominates $\delta$.

\eprop

\bproof Since normalization is a local operation, by Remark \ref{ru.r2.3} (2) we can view $X$ as a germ of  $\C^n$ at  the origin  with coordinates $(u, v, w_1 \ldots , w_{n-2})$
such that $D$ is given by $u=0$, $Z$ by $u=v=0$  (i.e. $\hX$ can be viewed as the hypersurface in $X \times \C_w$ given by the equation $uw=v$). 
By Lemma \ref{ref.l3} $g=ev^k +a$ where $a\in fI^{k-1}$ and $e$ is an invertible holomorphic function.
%That is, $f=u$, $h=v$, $\hX$ can be viewed as the smooth hypersurface given in $X \times \C_w$ by the equation $uw=v$. 
%In this case $fI^{k-1}$ is generated by functions $u^k, u^{k-1}v, \ldots , uv^{k-1}$. This implies the following.
%{\em Claim}. The Davis modification $\delta : \tX \to X$, given by the set-theoretical complete intersection $f^k=g=0$, is homogeneous if and only  $g$ is of form
%$g=ev^k +a+b$ where $a\in fI^{k-1}$, $b\in I^{k+1}$, and $e$ is an invertible holomorphic function. 
Changing $e$  we can suppose that the Taylor series of $a$ does not contain monomials divisible by $v^k$. 
Hence $u=g/f^k=g/u^k=ew^k +c$ where $c$ is a polynomial in $w$ of degree at most $k-1$ whose coefficients are holomorphic functions on $X$.
Thus $w$ is integral over the ring of holomorphic functions on $\tX$ which concludes the proof.

\eproof

\brem\label{ru.r3.5} (1) Note that for a simple pseudo-affine modification from Definition \ref{ru.d2a.3} (2) with a given divisor $D$ the function $f$ is determined uniquely
up to an invertible factor. 
%Furthermore, for a given center $Z$ to be a strict complete intersection $f=h=0$ the only allowed change of $h$
%is of form $eh+af$ where $e$ and $a$ are holomorphic functions on $X$ and $e$ is also invertible. 
This implies that for a given $k \geq 1$ the notion of  a homogeneous
modification of degree $k$ from  Definition \ref{ru.d2a.3} (2)  is determined  not as much by $f$ and $h$ but by $D$ and the defining ideal $I$ (or equivalently the center $Z$).

(2) In particular if $\delta : \tX \to X$ is a pseudo-affine modification with a divisor $D$ and a center $Z$ such that $\dim Z = \dim X -2$ then for smooth point $z$ of $Z$ that is not in $X_{\rm sing} \cup D_{\rm sing}$
we can say whether $\delta$ is locally homogeneous at $z$ or not (where in the former case there exists a Euclidean neighborhood $U$ of $z$ in $X$ such
that for every connected component $\tU$ of $\delta^{-1}(U)$ the modification $\delta|_\tU : \tU \to U$ is homogeneous).

\erem

\bprop\label{mod.new.2} Let $X$ be a normal irreducible  Stein space which is Cohen-Macaulay, $\delta : \tX \to X$ be a Davis modification with  an irreducible divisor
 $D$ and a center $Z$ of ${\rm codim}_X Z =2$ such that none of  irreducible component of $Z$ is contained in $X_{\rm sing} \cup D_{\rm sing}$.
 Suppose also that $Z$ is a strict complete intersection of form $f=h=0$ where $D=f^*(0)$ and that at general points of $Z$ this modification $\delta$ is locally homogeneous.
 Let $\sigma : \hX \to X$ be a simple pseudo-affine modification associated with divisor $D$ and center $Z$. Then $\hX$ is a normalization of $\tX$ and thus $\sigma$
 dominates $\delta$.
 
\eprop

\bproof First note that $\hX$ is a normal Stein space by Proposition \ref{an.normality}. 
Let $z$ be a general point of $Z$, i.e. $z$ is a smooth point of $Z$ and it is not contained in some subvariety $P\subset Z$ of ${\rm codim}_Z P \geq 1$ such that $Z\cap  (X_{\rm sing} \cup D_{\rm sing})\subset P$.
 %none of irreducible components of $Z$. 
Then by Proposition \ref{ru.p2.4} there exists a Euclidean
neighborhood $U$ of $z$ in $X$ such that $\sigma^{-1}(U)$ is a normalization of $\delta^{-1}(U)$. %Consider  a principal divisor $T$ in $X$ containing $P$ but not $z$.
%Hence   $\sigma^{-1}(X\setminus T)$ is a normalization of $\delta^{-1}(X \setminus T)$.
%In fact, varying $T$ we get 
That is, we have a biholomorphism $\psi$ between $\sigma^{-1}(X\setminus P)$ and a normalization of $\delta^{-1}(X \setminus P)$. % since $P$ is the intersection of all such $T$'s with $D$.
Since $P$ is of codimension at least 3 in $X$ and $\delta$ is a Davis modification the codimension of $\delta^{-1}(P)$  in $\tX$ is at least of  2 (because by Remark \ref{ru.r2.3}(3)  $\delta^{-1}(P)\simeq P \times \C$).
Hence the Hartogs theorem implies that $\psi$ extends to a biholomorphism between $\hX$ and a normalization of $\tX$ which is the desired conclusion. 

\eproof

\section{Semi-finite modifications}

\bnota\label{3.1}
Let $T$ be a germ of an analytic set at the origin $o$ of $\C^m$, $\Hol (T)$ be the ring of holomorphic functions on $T$, and
$D=T \times \C$. We consider a hypersurface $Z$ in $D$ such that  $\pi |_Z : Z \to T$ is finite where $\pi : D \to T$ is the natural projection.
We also suppose that $Z$ coincides (as a set) with the zeros
of an analytic function $h$.

\enota

\blem\label{ru.l1} Let Notation \ref{3.1} hold and  $w$  be a coordinate on the second factor of $D=T \times \C$. Then $h$ can be chosen as a monic polynomial in $w$ with coefficients from  $\Hol (T)$
(i.e.  $h \in \Hol (T)[w]$).

\elem

\bproof Let $\pi^{-1} (o) \cap Z$ consists of points $z_1, \ldots ,  z_n$ where  $z_i=(o,w_i) \in T\times \C$. By the Weierstrass preparation theorem for every $z_i$ there is
a Euclidean neighborhood in $D$ in which $h$ coincides with $e_i h_i(w)$ where $e_i$ is an invertible function and $h_i \in \Hol (T)[w]$ is a monic polynomial in $w-w_i$.
Note that $Z$ coincides with the zeros of the product $h_1 h_2 \cdots h_m$ because of the finiteness of   $\pi |_Z : Z \to T$.
Thus replacing $h$ with the product $h_1 h_2 \cdots h_m$ we get the desired conclusion.

\eproof

\blem\label{3.2} Let Notation \ref{3.1} hold and the zero multiplicity of $h$ at general points of $Z$ is $n$. 
Then $g=h^{1/n}$ is a holomorphic function on $D$ and in particular  the defining ideal of $Z$ in the ring $\Hol (T)[w]$ is the principal ideal 
generated by $g$.
\elem

\bproof  %Since  $\pi |_Z : Z \to T$ is finite, the leading coefficient of $h$ (as  a polynomial in $w$) is nowhere zero.
%Hence 
By Lemma \ref{ru.l1} we can suppose that $h$ is a monic polynomial in $w$. Consider the restriction of $h$ to any fiber. It is a polynomial
whose roots have multiplicities divisible by $n$ (since for general fibers such multiplicities are exactly $n$).
Thus  the restriction of $g$ to any fiber of $\pi : D \to T$ can be chosen as a nonzero monic polynomial in $w$. 
Therefore, choosing such a restriction on $\pi^{-1}(o)$ we define by continuity a unique branch of $h^{1/n}$ on $D$ as $g$. 
Note that $g$ is holomorphic outside a subset $K=Z \cap \pi^{-1} (T_{\rm sing})$
and because of the assumption on $h$ we see that  $g=w^k +r_{k-1}w^{k-1}+ \ldots + r_1w +r_0$ where each $r_i$ is
a continuous function on $T$ holomorphic on $T\setminus K$. 
Hence the first statement will follow from the claim.

{\em Claim}. Let $g$ be a (not  a priori continuous) function on $D$ holomorphic outside a proper analytic subset $K$
that does not contain any fiber of $\pi$. %and has codimension at least 2. 
Suppose also that $g$ is a polynomial in $w$. Then $g$ is holomorphic on $D$.

Taking a smaller $T$, if necessary, one
can suppose that the image $K_0$ of $K$ under the natural projection $D \to \C_w$ is relatively compact.
In particular, for every $w_0$ in $\C \setminus K_0$ the function $r_kw_0^k + r_{k-1}w_0^{k-1}+ \ldots + r_1w_0 +r_0$ is holomorphic on $T$.
Since we have infinite number of such $w_0$'s every coefficient $r_i$ is also holomorphic on $T$ 
which concludes the proof of the Claim.

To see that $Z$ is a principal divisor in $D$ consider any function $f\in \Hol (T) [w]$ vanishing on $Z$. Then the quotient $f/g$ is again holomorphic on $D$ by the Claim
which yields the desired conclusion. 

\eproof

%\brem\label{ru.r1} Note that $g$ is a the only monic polynomial in $w$ whose zeros coincides with the reduced divisor $Z$ even if we extend the ring of coefficients from $\Hol (T)$
%to the field of meromorphic functions on $T$.
%\erem

\blem\label{3.3} Let Notation \ref{3.1} hold and $T$ be singular. Then so is $Z$.

\elem

\bproof Assume that $Z$ is smooth at some point $z \in Z \cap \pi^{-1}(o)$. Let $g$ be a function as in Lemma \ref{3.2}.
Extend $g$ to a function $\tilde g \in \Hol (\C^m, o)[w]$ whose zeros define an extension $\tZ$ of $Z$. 
%By the Weierstrass preparation theorem the defining ideal $J$ of $\tZ$ at $z$ can be chosen as a principal ideal 
%generated by a function from $\Hol (\C^m, o)[w]$, i.e. we can suppose that this function is $\tilde g$ (taking a smaller 
%$Z$ if necessary).
Note that if the partial derivative of $\tilde g$ with respect to $w$ is nonzero at $z$ then by the implicit function theorem 
$\tZ$ is biholomorphic to $(\C^m, o)$ and therefore $Z$ is
locally biholomorphic to $T$ which implies that $Z$ is singular. Hence we assume that this derivative is zero. 
%Therefore, it is
%true for every function from $J$. 

Let $I$ be the defining ideal of $T$ in $\Hol (\C^m, o)$ and $k$ be the dimension of $T$. 
Choose any $m-k$ elements from $I$ and consider their Jacobi matrix
(with respect to coordinates $z_1, \ldots z_m$ of $\C^m$). 
Since $T$ is singular at $o$ any  $(m-k)$-minor $M$ of this matrix vanishes at $o$. 
%Hence the same remains true if we choose 
%these elements from $I [w]$.  

By Lemma \ref{3.2} the defining ideal $L$ of $Z$ in $(\C^m,o) \times \C_w$ is generated by $I$ and $\tilde g$.
Take any $m+1-k$ elements from $L$ and consider 
their Jacobi matrix with respect to coordinate $(w,z_1, \ldots z_m)$.
It follows from the previous argument about the partial derivative of $\tilde g$ with respect to $w$ and about
$(m-k)$-minor $M$ that any $(m+1-k)$-minor of this new Jacobi matrix vanishes at $o$.
This is contrary to the fact that $Z$ is smooth
at $z$ which concludes the proof.
\eproof

\bexa\label{ru.e3.1} Lemma \ref{3.3} is one of the central technical results in this paper. The assumption that $Z$ coincides with zeros of a global 
analytic function is very important. Indeed, consider the case when $T$ is a semi-cubic parabola $x^2-y^3=0$ in $\C^2_{x,y}$.
Then there is a closed immersion of $\C$ into $D=T \times \C \subset \C_{x,y,w}^3$ given by $t \to (t^3,t^2,t)$ such that the image is a smooth Weil divisor
whose projection to the singular $T$ is finite. Lemma \ref{3.3} is not applicable here because this Weil divisor is not $\Q$-Cartier.
\eexa

\brem\label{3.3a} If $T$ is not unibranch at $o$ then the argument is much easier and, furthermore, $Z$ is not unibranch at any point $z_0$ above $o$ as well.
Indeed, let $T_1$ and $T_2$ be distinct irreducible components of $T$ at $o$. Then $D$ has irreducible branches $T_1 \times \C$ and $T_2\times \C$
meeting along the line $L=o \times \C$. Since the zeros of $h$ contain the point $z_0 \in L$ but not the line $L$ itself the zeros $Z_i$ of $h$ in $T_i \times \C$
produce two different branches $Z_1$ and $Z_2$ of $Z$ at $z_0$.

\erem

\bprop\label{ru.p2}
Let  $T$ be an affine algebraic variety, $D=T \times \C$, and $Z$ be an algebraic hypersurface in $D$ such that
$\pi |_Z : Z \to T$ is finite where $\pi : D \to T$ is the natural projection. Suppose that 
for every point $t_0 \in T$
there is a Euclidian neighborhood $U\subset T$  such that $Z\cap \pi^{-1}(U)$ is an analytic $\Q$-principal divisor in $\pi^{-1}(U)$, i.e. for some natural $k>0$  and a holomorphic function $g$ on $\pi^{-1}(U)$
the divisor $kZ\cap \pi^{-1}(U)$ coincides with $g^*(0)$.

Then $Z$ is an effective reduced principal divisor in $D$.
\eprop

\bproof Since the field $\C (Z)$ of rational functions on $Z$ is a finite separable extension of the field $\C (T)$ of rational functions on $T$ there exists
$w \in \C (Z)$ such that $\C (Z)$ coincides with the field of rational functions in $w$ with coefficients from $\C (T)$. Let $g(w) = w^n + r_{n-1}(t) w^{n-1}+ \ldots + r_1(t) w+r_0(t)$
be the minimal monic polynomial for $w$ over $\C (T)$. 
In particular $Z \cap \pi^{-1}(U)$ is given by the zeros of this rational function $g$. By Lemma \ref{3.2} $Z \cap \pi^{-1}(U)$ is a principal
divisor in $\pi^{-1}(U)$ which implies that
the function $g$ is holomorphic. Therefore $g$ is regular (e.g., see \cite[Theorem 5]{Ka91}).
This yields the desired conclusion.

\eproof

\bdefi\label{3.4} Let $\sigma : \hX \to X$ be an affine modification of normal affine algebraic varieties such that  $X$ is Cohen-Macaulay and

(a)  its divisor $D=f^*(0)$ (where  $f \in \C [X]$) is a reduced principal divisor in $X$ isomorphic to a direct product
$D \simeq T \times \C$;

(b) the restriction $\pi|_Z : Z  \to T$ of the natural projection $\pi : D \to T$ to the center $Z$ of $\sigma$ is finite;

(c) none of irreducible components of $Z$ is contained in the singularities of $X$ (note that for the singularities of $D$ the similar fact also holds);

(d) for any point $z$ of $Z$ there is a Euclidean neighborhood $U$ in $X$ for which the restriction $\sigma|_\hU : \hU \to U$  of $\sigma$ to $\hU=\sigma^{-1} (U)$ 
factors through some Davis  pseudo-affine modification $\delta : \tU \to U$ with the following property: the restriction of $\delta$ over a neighborhood $U'\subset U$ of any general point 
$z' \in Z \cap U$ is homogeneous of degree $k$ (where $k$ does not depend on different irreducible components of $Z \cap U$).

Then we call such a $\sigma$ a semi-finite  affine modification.

\edefi

\bprop\label{3.5} Let $\sigma : \hX \to X$ be a semi-finite affine modification with  $D \simeq T \times \C$ and $Z$ as in Definition \ref{3.4}.

Then $\sigma$ dominates a simple modification (with the same center and divisor) and if $T$ is singular so is the center $Z$.
\eprop

\bproof  Let a Davis modification $\delta : \tU \to U$  from Definition \ref{3.4} (d) be defined by  the ideal $(f^{k},g)$. By Condition (d)  $\delta$ is locally homogeneous at general points of $Z\cap U$
which in combination with Lemma \ref{ref.l3}  % the Claim in Proposition \ref{ru.p2.4}   
implies that for some $k\geq 1$ the divisor $kZ\cap U$ coincides with $g^*(0)\cap U$. 
That is,
 we have Condition (i) from Proposition \ref{ru.p2}.
Hence  $Z=h^*(0)$ where $h\in \C [D]$.
Thus extending $h$ to a regular function on $X$ (denoted by the same symbol) we see that $Z$ is a strict complete intersection given
by $f=h=0$. These functions $f$ and $h$ induce a simple affine modification $\tau : X' \to X$ with divisor $D$ and center $Z$
such that $X'$ is a normal affine algebraic variety by Proposition \ref{normality}. Furthermore, $\hX$ and $X'$ are also normal as analytic sets by \cite[Chapter 13, Theorem 32]{ZaSa}.
By Proposition \ref{mod.new.2} 
$\tau^{-1}(U)$ is an  analytic normalization of $\tU$. Since $\hX$ is normal then the holomorphic map $\sigma^{-1}(U) \to \tU$ (induced by the domination of $\delta$ by $\sigma$)
factors through the normalization of $\tU$. This yields a desired holomorphic map from $\hX$ to $X'$ which is a morphism since its restriction over $X\setminus D$
is an algebraic isomorphism.  That is, $\sigma$ dominates $\tau$. The second statement follows from
Lemma \ref{3.3}.

\eproof

\section{Applications of Kambayashi's theorem}

\bnota\label{ru.n5.1} In this section $\varphi : X \to B$ will be a morphism of complex factorial affine algebraic varieties,
$\Phi$ will be a $G_a$-action $\Phi$ on $X$ which preserves the fibers of $\phi$. We do not assume a priori that the $\C [B]$-algebra
$\C [X]^\Phi$ of $\Phi$-invariant regular functions is finitely generated.
\enota

However, such an algebra can be viewed  (and will be viewed) as a direct limit $\varinjlim A_\alpha$ of its finitely generated $\C [B]$-subalgebras $A_\alpha$ (with respect to the partial order
generated by inclusions) where $\alpha$ belongs to some index set and each
$A_\alpha$ can be treated as the ring $\C [Q_\alpha]$ of regular functions of some affine algebraic variety $Q_\alpha$.
Replacing $A_\alpha$ with its integral closure we suppose that each of these $Q_\alpha$ is normal
(the fact that this transition preserves affineness is a standard result, e.g. \cite[Theorem 4.14]{Ei}).
%Inclusions $A_\alpha \subset A_\beta$ generates morphisms $Q_\beta \to Q_\alpha$ which make
%${\rm Spec} \, \C [X]^\Phi$ the inverse limit of such $Q_\alpha$'s. 
Furthermore, by the Rosenlicht Theorem (e.g., see \cite[Theorem
2.3]{PV}) $Q_\alpha$ can be chosen so that the morphism $\rho_\alpha : X \to Q_\alpha$ induced by the natural embedding
$A_\alpha \subset \C [X]$ separates general orbits of $\Phi$. As in \cite{GW} we introduce the following. 

\bdefi\label{ru.d5.1}
A normal affine variety $Q_\alpha$ as before will be called a
{\em partial quotient} of $X$ by  $\Phi$ and the morphism $\rho_\alpha : X \to Q_\alpha$ (separating general orbits of $\Phi$) will be called {\em a partial quotient morphism}. \edefi

\brem\label{ru.r5.1} Note that including the coordinate functions of morphism $\varphi$ into a ring $A_\alpha$ we can always choose a partial quotient morphism that
is constant on the fibers of $\varphi$. 

\erem

\bthm\label{4.1} Let $\varphi : X \to B$ and $\Phi$ be as in Notation \ref{ru.n5.1} and  the generic fiber of $\varphi$ be a three-dimensional variety $Y$ (over the field $K$
of rational functions on $B$).  Let $\tK$ be the algebraic closure $K$
and $\tY$ be the variety over $\tK$ obtained from $Y$ by the field extension. Suppose that
$\rho : X \to Q$ is a partial quotient morphism such that $\varphi$ factors through $\rho$. Suppose also
that  the ring of invariants of
the $G_a$-action on $\tY$ (induced by $\Phi$) is a polynomial ring in two variables over $\tK$.  Then there is a
morphism $\psi : Q \to B\times \C^2$ such that for some nonempty Zariski dense open subset $B^* \subset B$
the restriction of $\psi$ over $B^*$ is an isomorphism. 

\ethm

\bproof By the Kambayashi theorem \cite{Kam} the ring of invariants of the induced action on $Y$ is 
also a polynomial ring in two variables over $K$. Hence for some $B^*$ as above and
$X^* = \varphi^{-1} (B^*)$ the induced action on $X^*$ has the categorical quotient isomorphic to $B^* \times \C^2$.
%By the universal property of quotient morphisms there is a map $\iota : Q^* \to Q$.
Assume that $B \setminus B^*$ is
a principal effective divisor which can be done without loss of generality.
Let $f$ be a regular function on $B$ with zero locus $B \setminus B^*$ (we denote the lifts of $f$ to $X$ or $Q$ by the same symbol). For every regular function $h \in \C [X^*]$
there exists  natural $k$ for which $f^kh$ extends to a regular function on $X$.  Hence for sufficiently large $k_1$ and $k_2$
the composition of the quotient morphism $\rho_0: X^* \to  B^* \times \C^2_{u_1,u_2}$ with the isomorphism  $\kappa : B^* \times \C^2_{u_1,u_2}\to  B^* \times \C^2_{u_1,u_2}$ 
given by $(b,u_1,u_2) \to (b,f^{k_1}u_1, f^{k_2}u_2)$ extends to a morphism $\tau : X \to B \times \C^2$ (note that $\tau|_{X^*} : X^* \to B^* \times \C^2$ can be
viewed now as the quotient morphism).
Similarly, since $f\in \C [X]$ and each $u_i\in \C [X^*]$  are $\Phi$-invariant we see 
that $f^{k_1}u_1$ and $f^{k_2}u_2$ can be viewed as regular functions on the normal variety $Q$.
Thus
$\tau =\psi \circ \theta $ where $\theta : X \to Q$ and  $\psi : Q \to B\times \C$ are morphisms. %while the restriction of $\tau$ to $X^*$ coincides with the quotient morphism $\kappa \circ \rho_0 : X^* \to B^* \times \C^2$.
By the universal property of quotient morphisms $\theta |_{X^*} : X^* \to Q \setminus f^{-1}(0)$ factors through $\tau |_{X^*}$ which implies that $\psi$ is invertible over $B^*$.
This yields the desired conclusion.
\eproof

\brem\label{4.2} By construction $\psi$ is an affine modification. Another observation that if the general fibers of $\varphi$ do not
admit nonconstant invertible functions then for every irreducible divisor $T$ in $B$ the variety $\varphi^{-1} (T)$ is reduced and
irreducible by Proposition \ref{fv.p1}. Actually, in this case there is no need to assume that $B$ is factorial
since it follows automatically from the fact that $X$ is factorial.
\erem

\bexa\label{4.3}  Let Notation \ref{ru.n5.1} hold and
every general fiber of $\varphi$ be isomorphic to the same affine algebraic threefold $V$. Suppose that 
$\Phi$ induces an action on each
general fiber of $\varphi$  with the categorical quotient isomorphic to $\C^2$.
%\footnote{In the original version of the paper there was an additional assumption
%that general fibers of $\varphi$ are pairwise isomorphic, but the referee observed that this assumption is superfluous.}
Then we claim that
the above conclusion about the partial quotient
$Q$ as a modification of $B \times \C^2$ over $B$ is valid.

Indeed, by \cite{KrRu}  there exists a Zariski dense open subset $B^* \subset B$ and an unramified
covering $\hB^* \to B^*$ such that $\hX^*=X \times_{B^*} \hB^*$ is naturally isomorphic to $\hB^* \times V$ over $\hB^*$.
In particular, the induced action on the generic fiber of $\hX^* \to \hB^*$ has the ring of invariants isomorphic
to the polynomial ring $\hK [x_1, x_2]$ where $\hK$ is the field of rational functions on $\hB^*$.
Note that this ring is obtained from the ring of invariants on the generic fiber $Y$ of $X^* \to B^*$ via a field
extension $[\hK : K]$. Hence  by the Kambayashi theorem the latter ring of invariants is $K [x_1,x_2]$ and
we under the assumption of Theorem \ref{4.1}.
\eexa

In fact, as it was suggested by the referee, we can strengthen Theorem \ref{4.1} and Example \ref{4.3} as follows.
%In fact the same combination of the Kraft-Russell theorem and the Kambayashi theorem yields a stronger result and we have
%the following theorem which belongs to the referee.

\bthm\label{ref.t4}  Let Notation \ref{ru.n5.1} hold and let
$\Phi$ induce an action on each
general fiber of $\varphi$  with the categorical quotient isomorphic to $\C^2$. Let $Y$ and $K$ be as in Theorem \ref{4.1}.
Then the ring of invariants of the $G_a$-action on $Y$ (induced by $\Phi$) is a polynomial ring in two variables over $K$.
In particular, for a partial quotient morphism $\rho : X \to Q$ for which $\varphi$ factors through $\rho$ (i.e. for some morphism $\kappa : Q \to B$ one has $\varphi = \kappa \circ \rho$) 
there is a morphism $\psi : Q \to B\times \C^2$  over $B$ such that for a nonempty Zariski dense open subset $B^* \subset B$
the restriction of $\psi$ over $B^*$ is an isomorphism. \ethm

In preparation for the proof of this theorem we need the next (certainly well-known) fact.

\bprop\label{ref.p41} Let $\kappa : Q \to B$ be a dominant morphism of algebraic varieties such that $Q$ is normal.
Then for a general point $b \in B$ the variety $\kappa^{-1}(b)$ is normal.
\eprop

\bproof Replacing $B$ by a Zariski dense open subset (and $Q$ by its preimage) we can suppose that $\kappa$ is flat  \cite[Th\'eor\`eme 6.9.1]{EGA}.
Let $\omega$ be a generic point of $B$. Then the fiber of $\kappa$ over $\omega$ is normal (e.g., see \cite[Proposition 5.13]{AM})
which is equivalent to the fact that this fiber is geometrically normal\footnote{Recall that a scheme $Z$ over a field $k$ is geometrically normal if it is normal and, furthermore, it remains normal under any field extension of
$k$. In the case of a perfect field $k$ (e.g. a field of characteristic zero) every normal scheme is automatically geometrically normal.}  since we work over a field of characteristic zero. On the other hand the
set of  (not necessarily closed) points in $B$ for which the fibers over them are geometrically normal is open  \cite[Th\'eor\`eme 12.1.6]{EGA3}.
Since this set is nonempty it contains general (closed) points of $B$ and we are done.
\eproof

\subsection{Proof of Theorem \ref{ref.t4}.}
Let $Q_b$ be the fiber 
 $\kappa^{-1}(b)$ over  a general point $b\in B$ and $Q^b \simeq \C^2$ be the categorical quotient of the action $\Phi|_{\varphi^{-1}(b)}$. 
By the universal property of quotient
morphisms one has a morphism $\psi_b : Q^b \to Q_b$. Since $\rho$ separates general orbits $\psi_b$ must be birational. 

Assume that for a curve $C_b \subset Q^b$
the image $\psi_b(C_b)$ is a point $q_b \in Q_b$. By Corollary \ref{fv.c1} the preimage of $C_b$ in $\varphi^{-1}(b)$ is a surface $S_b$, i.e. $\rho (S_b) =q_b$.
By  \cite[Chapter I, Section 6.3, Corollary]{Sha} the closure  $T$ of the set $\{ q \in Q | \dim \rho^{-1}(q)=2 \}$ is a subvariety. Let $L$ be the union of the one-dimensional components of $T$
(it is non-empty  since $q_b\in T$  for general $b$).
Note that $\rho^{-1}(L)$ is a closed $\Phi$-stable threefold $V$.  
%By Proposition \ref{ref.p1} $f$ is a lift of a function $g \in \C [Q]$. 
By  the Hironaka flattening theorem \cite{Hi} there exists  of a proper birational morphism $\pi : \hQ \to Q$
such that for the irreducible component $\hX$ of  $X \times_Q \hQ$ dominant over $X$ the natural projection $\hat \rho : \hX \to \hQ$ is flat, i.e. all of its fibers are one-dimensional or empty. 
This implies $\pi^{-1}(L)$ contains
a surface $\hP =\hat \rho (\hV)$ where $\hV$ is the proper transform of $V$ in $\hX$.
Choose a curve $\hC\subset \hP$ whose image in $B$ is dense and such that $\hat \rho^{-1}(\hat c) \ne \emptyset$ for a
general point $\hat c \in \hC$.  Let $\hR$ be  a closed surface in $\hQ$ that meets  $\hP$ along $\hC$ and let $R$ be the closure of $\pi (\hR )$  (note that $R$ contains $L$).
Then the closure of $\rho^{-1} (R\setminus L)$ contains a $\Phi$-stable threefold $W$ that meets $V$ over general points of $B$ (since $\hat \rho^{-1}(\hat c) \ne \emptyset$). 
Because $X$ is factorial the threefold
$W$ is the zero locus of a regular function $h \in \C [X]$ which, by construction, does not vanish on general fibers of $\rho$. That is,  it is constant on general fibers and, therefore, $\Phi$-invariant.  
By construction $h$ is not constant on the surface  $V\cap \varphi^{-1}(b)$ for general $b$ (this follows from the fact that $W$ meets $V\cap \varphi^{-1}(b)$ along a curve).
Enlarging a set of generators of the ring $\C [Q]$ by $h$ (and taking the integral closure to preserve the normality of $Q$ assumed in Definition \ref{ru.d5.1})
we obtain another partial quotient for which the image of $V$ is a surface.  That is, the set $T$ as before becomes at most finite and for a general $b \in B$
this curve $C_b$ does not exists.

Hence the morphism $\psi_b : Q^b \to Q_b$ is quasi-finite now in addition to being birational. Assume that it is not an embedding. Then by the Zariski Main theorem $Q_b$ is not normal for general $b \in B$
contrary to Proposition \ref{ref.p41}.
%Let $a_b$ be a singularity of $Q_b$. 
%Since general points of $B$ are regular values of morphism $\kappa$ there is a vector $v\in T_{a_b}Q$ such that $\mu=\kappa_* (v)$ is a
%nonzero vector from $T_bB$. Choose a vector field on $Q$ whose value at $a_b$ is $v$. Dividing this vector field by a holomorphic function in a neighborhood $O\subset Q$ of $a_b$ we can obtain
%a vector field $\nu$ on $O$ with the following properties: $\nu (a_b)=v$ and for every $a'$ and $a'' \in O$ with $\kappa (a')=\kappa (a'')$ one has $\kappa_* (\nu (a'))=\kappa_* (\nu (a''))$.
%Applying the phase flow of this field one can see that a neighborhood of $a_b$ in $Q$ is biholomorphic to $O_b \times \Delta$ where $O_b$ is a neighborhood of $a_b$ in $Q_b$
%and $\Delta$ is the unit disc $\{ z \in \C | |z| <1 \}$. Hence $Q$ is not normal in the holomorphic sense and, therefore, it is not normal in the algebraic setting \cite[Chapter 13, Theorem 32]{ZaSa}. 
%This contradicts normality of $Q$.
Hence $\psi_b : Q^b \to Q_b$ is an embedding. 

Let $D_b$ be the complement to the image of $Q^b$ in $Q_b$. Suppose that $D_b$ is a curve for a general $b \in B$ (since $Q_b$ and $Q^b$ are affine the alternative is an empty $D_b$). 
The set $Q_b\setminus \rho (\varphi^{-1}(b))$ consists of $D_b$ and a finite set by Corollary \ref{fv.c1}.
Since $\rho (X)$ is a constructible set by \cite[Chap. II, Exercise 3.19]{Har}
we see that $Q \setminus \rho (X)$ is an algebraic variety. Consider the irreducible components of this variety which are dominant over $B$ and remove those of them that have dimension $\dim B$.
Then we are left  with $D:=\bigcup_{b \in B} D_b$
because  $Q_b\setminus (\rho (\varphi^{-1}(b)) \cup D_b)$ is finite. That is, $D$ is an algebraic variety and $Q\setminus D$ is a quasi-affine variety. 
The general fibers $Q_b \setminus D_b$ of the natural morphism $Q\setminus D\to B$ are isomorphic to $\C^2$. By \cite{KaZa01} for a Zariski dense open subset $B^*\subset B$ the variety
$Q^*=\kappa^{-1} (B^*)\setminus D$ naturally isomorphic to $B^* \times \C^2$. Hence for $X^*=\varphi^{-1}(B^*)$ we get a partial quotient morphism 
$\rho|_{X^*} : X^* \to B^* \times \C^2$. 

Let $\rho' : X^* \to Q'$ be another partial quotient morphism for $\Phi|_{X^*}$  into a normal variety $Q'$ (over $B^*$)
such that $\rho$ factors through $\rho'$, i.e. $\rho = \theta \circ \rho'$ for a morphism $\theta : Q' \to B^* \times \C^2$. Suppose that $Q_b'$ is the fiber of
the natural morphism $Q' \to B^*$ over a point $b \in B^*$, i.e. the restriction of $\theta$ yields a morphism $\theta_b : Q_b'\to Q_b\simeq Q^b\simeq \C^2$.
By the universal property of quotient morphisms we have a natural morphism $Q^b \to Q_b'$, i.e.  $\theta_b$ is invertible. Hence
$\theta$ is bijective. By the Zariski Main theorem $\theta$ is an isomorphism.  This implies that $\rho|_{X^*} : X^* \to B^* \times \C^2$ is the categorical quotient morphism.
In particular, for $Y$ and $K$ as in Theorem \ref{4.1} the ring of invariant of the $G_a$-action on $Y$ (induced by $\Phi$) is a polynomial ring in two variables over $K$.
Now the desired conclusion follows from Theorem \ref{4.1}. \hspace{10.2cm} $\square$

\section{Criterion for  existence of an affine quotient}

\bnota\label{5.1} % In this section Notation \ref{ru.n5.1} holds with 
Let $B$ be a unibranch germ of a smooth complex algebraic curve at point $o$ (i.e. $o$ is
the zero locus of some $f \in \C [B]$) and $\varphi : X \to B$ be a morphism from a complex factorial affine algebraic variety $X$ equipped with a $G_a$-action $\Phi$
which preserves each fiber of $\varphi$. 
We  suppose also that $\rho : X\to Q$ is a partial quotient morphism of the action $\Phi$. % which preserves every fiber of $\varphi$. 
Let $\psi : Q\to Q_0$ be a morphism
over $B$ into a smooth affine algebraic variety $Q_0$ over $B$ and $\tau= \psi \circ  \rho$. The divisor in $X$ (resp. $Q$, resp. $Q_0$)
over $o$ will be denoted by $E$ (resp. $D$, resp. $D_0$). We suppose
that the morphism $Q_0 \to B$ is smooth
(and in particular $D_0$ is a smooth reduced divisor) and
that $\psi$ induces an isomorphism $Q \setminus D \to Q_0 \setminus D_0$. We suppose also that $E=\varphi^*(o)$ is reduced and denote by
$Z_0$ the closure of $\tau (E)$ in $D_0$. 
%and that every component of the closure $Z_0$ of $\tau (E)$ is of codimension 1 in $D_0$.
We consider the Stein factorization (e.g., see \cite[Chapter III, Corollary 11.5]{Har}) of a proper extension $\bE \to Z_0$ of the morphism $\tau|_E: E \to Z_0$.
Its restriction to $E$ enables us to treat $\tau|_E: E \to Z_0$ as a composition of  a surjective morphism
$\lambda :E \to V$ with connected general fibers and a  quasi-finite morphism $\theta : V\to Z_0$.  

\enota

Theorem \ref{1.1} implies the following.

\blem\label{5.2}  Let Notation \ref{5.1} hold and $Z_0$ (and therefore $\psi (D)$) be Zariski dense in $D_0$. Then $\psi$ is an isomorphism.

\elem

\bconv\label{5.2a} With an exception of Corollary \ref{5.7}  we suppose throughout this section that $Z_0$ is a divisor in $D_0$ and furthermore

(i) the morphism $\theta: V \to Z_0$ from Notation \ref{5.1} is in fact {\bf finite};

(ii) for every point $z_0\in Z_0$ the preimage $\tau^{-1} (z_0)$ is a {\bf surface}.

%(iii)  $E$ is smooth \footnote{This assumption can be replaced by the following:
%normalization $E_{\rm norm}$ of $E$ is smooth and for every holomorphic function $g\in \Hol (E_{\rm norm})$ there is $k>0$ for which 
%$g^k\in \Hol (E)$.}.

\econv

\brem\label{5.2b} Note that Convention \ref{5.2a} holds automatically when $E$ is isomorphic to $\C^3$ (under the assumption that
$Z_0$ is a curve). Indeed, $\tau^{-1}(z_0)$ cannot be three-dimensional (otherwise it coincides with $E$ and $\tau (E)$ is a point, not a curve).
Then both $V$ and  $Z_0$ must be polynomial curves and a nonconstant morphism of polynomial curve is always finite.
In fact, this argument works not only when $E\simeq \C^3$. It is enough to assume that there is no nonconstant morphism from $E$ into
any non-polynomial curve.
\erem

\bnota\label{5.3}

Let Notation \ref{5.1} hold. Our aim is to present $\psi : Q \to Q_0$ over $B$ as a composition

\bdi Q=:Q_n & \rTo^{\sigma_n}& Q_{n-1}&\rTo^{\sigma_{n-1}}& \ldots &\rTo^{\sigma_2} & Q_1 &\rTo^{\sigma_1}& Q_0
\edi

\noindent of simple affine modifications $\sigma_i : Q_i \to Q_{i-1}$ over $B$. Let 
$\psi_i = \sigma_1 \circ \ldots \circ \sigma_i : Q_i \to Q_0$ and suppose that $\psi$ can be presented as $\psi = \psi_i \circ \delta_i$ for a modification  $\delta_i : Q \to Q_i$.
Let $\tau_i = \delta_i \circ \rho$ and thus $\tau = \psi_i \circ \tau_i$. We suppose that the zero locus of $f\circ \psi_i$ is reduced 
and coincides with $D_i:=\psi_i^{-1} (D_0)$ (which is nothing but the divisor of modification $\sigma_{i+1}$) while the center $Z_i$ of $\sigma_{i+1}$
coincides with $\tau_i (E)$. Note that  under such assumption every $Q_i$ is normal and Cohen-Macaulay by Proposition \ref{normality} and induction starting
from smooth $Q_0$ and $D_0$.

\enota 

\blem\label{5.4} Suppose that for some $i$ one has a sequence of simple affine modifications
\bdi Q_i & \rTo^{\sigma_i}& Q_{i-1}&\rTo^{\sigma_{i-1}}& \ldots &\rTo^{\sigma_2} & Q_1 &\rTo^{\sigma_1}& Q_0 \edi
\noindent such that $\psi : Q \to Q_0$ factors through $\psi_i : Q_i \to Q_0$ and, therefore, $\psi = \psi_i \circ \delta_i$. Let $\delta_i : Q \to Q_i$ be a semi-finite
affine modification. 
%and that there exists a semi-finite Davis
%modification $\tilde \sigma_{i+1} : \tQ_{i+1} \to Q_i$ such that $\delta_{i} : Q\to Q_i$ factors through $\tilde \sigma_i$.
Then the sequence can be extended to a sequence of simple affine modifications 
\bdi Q_{i+1} & \rTo^{\sigma_{i+1}}& Q_{i}&\rTo^{\sigma_{i}}& \ldots &\rTo^{\sigma_2} & Q_1 &\rTo^{\sigma_1}& Q_0 \edi
\noindent such that $D_i$ (resp. $Z_i$) is the divisor (resp. center) of $\sigma_{i+1}$ and $\psi : Q \to Q_0$ factors through $\psi_{i+1} : Q_{i+1} \to Q_0$, i.e. $\psi=\psi_{i+1} \circ \delta_{i+1}$.

\elem

\bproof The statement follows from Proposition \ref{3.5}. 

\eproof

In order to establish when  $\delta_i$ is semi-finite we need the following technical fact.

\blem\label{technical} For any point $x_0 \in E$ there is the analytic germ $S$ at $x_0$ of an algebraic subvariety of $X$ of dimension $\dim X -2$ such that

{\rm (i)} $S_0:=\tau (S)$ is a germ of an analytic  hypersurface in $Q_0$ at $z_0=\tau (x_0)$;

{\rm (ii)} the map $\tau|_S : S \to S_0$ is finite; 

{\rm (iii)} for every general point $z \in Z_0$ near $z_0$ there is a local coordinate system $(u,v,\bar w)$ on $Q_0$ at $z$ for which
$D_0$ is given locally by $u=0$, $Z_0$ by $u=v=0$, and every irreducible branch of the germ of $S_0$ at $z$ by an equation of form $e:=v-u^ld(u, \bar w)=0$ where $d$ is holomorphic
and $l$ is the zero multiplicity of the function $e \circ \tau$ at general points of $E$.

\elem

\bproof 

Consider two regular functions $h$ and $e$ on $E$ with the set $P$ of common zeros
such that near any general point $P$ is locally a strict complete intersection given by $h=e=0$
and for $z_0=\tau (x_0)$ the set $P \cap \tau^{-1} (z_0)$ contains $x_0$ as an isolated component (which is possible since $\tau^{-1}(z_0)$ is
a surface).
Extend $h$ and $e$ regularly to $X$. Without loss of generality we can suppose that
the set of common zeros of these extensions is an $(n-2)$-dimensional subvariety $T$ of $X$ where $n=\dim X$. Denote by $S$ (resp. $R$) the analytic germ
of $T$ (resp. $P$) at $x_0$ and by $T_0$ the closure of $\tau (T)$ in $Q$. Let $S_0$ be  the analytic germ of the  hypersurface $T_0 \subset Q_0$ at $z_0$.
For the restriction  $\kappa : S\to T_0$ of $\tau$
the preimage $\kappa^{-1}(z_0)=x_0$ is a singleton.  Hence \cite[Chapter 1, Sec. 3, Theorem 2]{GrRe} implies that the holomorphic map $\tau|_S : S\to S_0$ is finite.
Thus we have (i) and (ii).

%by the semi-continuity theorem \cite[Chapter I, Section 6.3, Corollary]{Sha} $\kappa$ is quasi-finite.
%We can suppose that
%the germs $R\subset S$ are relatively compact in $X$. By construction, the image $\tau (\partial R)$ of the boundary $\partial R$ of $R$  does not contain $z_0$ and, therefore, 
%$\tau (\partial S)$  does not contain $z_0$ as well (where  $\partial S$ is the boundary of $S$).
%Choose a small relatively compact neighborhood $S_0$ of $z_0$ in $T_0$ such that
%$\bar S_0$ does not meet $\tau (\partial S)$. Then the preimage
%$\kappa^{-1} (\bar S_0)$ in $S$ is compact.  Replacing $S$ with $\kappa^{-1} (S_0)$ we see that 
%$\tau|_S : S\to S_0$ is proper. By Grothendieck's theorem \cite[EGA IV, 8.11.1]{EGA} in combination with quasi-finiteness this implies that $\tau|_S : S\to S_0$ is finite.
%This yields (i) and (ii).

By construction  $S$ meets $E$ transversely at  a general point $x$ of $R$.
Let $z=\tau (x)$ and $(u',v',v'',\bar w')$ (resp. $(u,v, \bar w)$) be a local analytic coordinate system at $x \in X$ (resp. $z\in Q_0$) as
in Lemma \ref{coordinates}. That is,  locally $E$ (resp. $D_0$, resp. $Z_0$) is given by $u'=0$ (resp. $u=0$, resp. $u=v=0$) and $ \tau$ is given by $u=u'$, $\bar w = \bar w'$,
and $v=(u')^lq(u',v',v'', \bar w')$. 

The finiteness of $\tau|_S$ implies that $\tau|_{S\cap E} : S\cap E \to S_0\cap D_0$ is  \'etale over the general point $z \in Z_0$.
Thus $E\cap S$ is given locally by equations of form $v'=h_1(\bar w') $ 
and  $v''=h_2(\bar w')$ where $h_1$ and $h_2$ are holomorphic functions.
Since $S$ meets $E$ transversely we also see that $S$ must be given locally by equations of form $v'=h_1(\bar w') + u' g_1(u',v',v'', \bar w')$ 
and  $v''=h_2(\bar w') + u' g_2(u',v',v'', \bar w')$ where $g_1$ and $g_2$ are holomorphic functions. Furthermore, by the implicit function theorem these equations
can be rewritten as $v'=h_1(\bar w') + u'\tilde g_1(u', \bar w')$ 
and  $v''=h_2(\bar w') + u'\tilde g_2(u', \bar w')$ for some other holomorphic functions $\tilde g_1$ and $\tilde g_2$. Plugging these expressions for $v'$ and $v''$
with $u'=u$ and $\bar w=\bar w'$ 
into equation $v=(u')^lq(u',v',v'', \bar w')$ we get the desired form of equation $e=0$ in (iii). 
Note that the function $e \circ \tau$  is given by $$e \circ \tau (u',v',v'',\bar w') =(u')^l [q(u',v',v'',\bar w')-h(u', \bar w')].$$
Since $q(0,v',v'',\bar w')$ depends on $v'$ or $v''$ by Lemma \ref{coordinates},
we see that the expression on in the brackets is not identically zero on $E$ which shows that $l$ is the zero multiplicity of $e \circ \tau$ at a general point of $E$.

\eproof

%\brem\label{ru.r5.5} The only restrictions on the choice of, say, function $h$ in the proof of Lemma \ref{technical} are that  $W=h^{-1}(0)$ contains $x_0$, that $h$ does not vanish identically
%on any component of $\tau^{-1}(z_0)$, and that the divisor $h^*(0)$ is reduced. 
%%at general points of $V=h^{-1}(0)$ it has  simple zeros and $x_0 \in h^{-1}(0)$.
%Then the second desired function $e$ can be always constructed.
% Thus, fixing such $h$ we can always suppose that
%$S$ is contained in a given principal hypersurface $W$.
%\erem

\blem\label{5.5} The morphism $\delta_i$ is always semi-finite 
%A semi-finite Davis modification $\tilde \sigma_{i+1}: \tQ_{i+1} \to Q_i$ as in Lemma \ref{5.4} always exists 
unless $Z_{i}=\tau_i(E)$ is dense in $D_i$.
\elem

\bproof As we mentioned in Notation \ref{5.3} the variety $Q_i$ is normal and Cohen-Macaulay, and the divisor $D_i=(f\circ \psi_i)^*(0)$ is principal.
Furthermore, by construction $D_i=Z_{i-1}\times \C$ for $i\geq 1$. Thus we have condition (a) from Definition \ref{3.4}.
The finite morphism $V \to Z_0$  (from Notation \ref{5.1} and Convention \ref{5.2a}) factors through maps
$\theta_i : V \to  Z_{i} $, $Z_i\to Z_{i-1}$ and $Z_{i-1} \to Z_0$ each of which must be, therefore, finite, i.e. we  have condition (b) from Definition \ref{3.4}.
Condition (c) follows from construction and Proposition \ref{normality}. 
It remains to check conditions (d). %where without loss of generality we suppose that $i=0$ in Lemma \ref{5.4}, i.e. $Q_i=Q_0$.
%Since this condition is local we can replace $Q_0$ by a Euclidean
%neighborhood of a given point $z_0 \in Z_0$ in $Q_0$ and consider the restrictions of $\psi, \psi_i$, $\delta_i$, and $\tau_i$ over it. 

Consider the following construction of an analytic germ of an algebraic variety $\cS$ in $X$ of $\dim \cS =n-2$ where $\dim X=n$. 
Let $M$ be the set of points in $V$ above $z_0$ (by Convention \ref{5.2a} $M$ is finite and nonempty).
For every $y \in M$ choose any $x\in E$ above $y$. Consider the analytic germ $S(y)$ of an $(n-2)$-dimensional algebraic subvariety of $X$ at $x$ such that $S_0(y):=\tau (S(y))$ is an analytic hypersurface in $Q_0$
and the restriction of $\tau$ yields a finite map $S(y) \to S_0(y)$ (the existence of such an $S(y)$ is provided by Lemma \ref{technical}). 

Let $V(y)$ be the analytic germ of $V$ at $y$ and $k(y)$ be the degree of the finite morphism $E\cap S(y) \to V(y)$ induced by $\lambda : E \to V$ from Notation \ref{5.1}.
In general $k(y)$ may depend on $y$ but we allow $S(y)$ to be non-reduced and then replacing each $S(y)$  (and, therefore, $S_0(y)$)
with a multiple of it we can suppose that $k(y)=k$ for every $y \in M$.\footnote{We use here the following definition: if  $r$ is the degree of a finite morphism $W \to U$ of reduced analytic sets
and $W_m$ is the $m$-multiple of $W$ then we put the degree of the morphism $W_m \to U$ equal to $mr$. Treating $W_m$ as a union of $m$ disjoint samples
of $W$ one can see that this extended notion of degree has all the properties of the standard degree.} Since the variety $Q_0$ is smooth every $S_0(y)$ coincides with $h_{0y}^*(0)$ for
some holomorphic function $h_{0y}$.  By Lemma \ref{technical} for a local analytic coordinate system $(u,v,\bar w)$ at a general point $z \in Z_0$
this function $h_{0y}$ can be viewed as $(v-u^ld(u,\bar w))^k$ where $l$ is the zero multiplicity of $(v-u^ld)\circ \tau$ on $E$.

Let $\cS =\bigcup_{y\in M} S(y)$, $\cS_0=\bigcup_{y\in M} S_0(y)$, and  $h_0=\prod_{y\in M} h_{0y}$, i.e. $\cS_0=h_0^*(0)$.
Then by construction  $\tau |_{\cS} : \cS \to \cS_0$ is finite and $\cS_0$ meets $D_0$ along the analytic germ of $Z_0$ at $z_0$.
For $m$ being the degree of the finite morphism $\theta : V \to Z_0$ from Convention \ref{5.2a} this implies that near point $z$ the function 
$h_0$ can be viewed as a product of $km$ factors of form $v-u^ld(u,\bar w)$. Hence the zero multiplicity of $h_0 \circ \tau$ at general points of $E$ is $klm$.

Let $S_i(y) =\tau_i(S(y))$ and $\cS_i=\tau_i (\cS)=\bigcup_{y\in M} S_i(y)$.  Note that 
the finiteness of  $\tau |_{\cS} : \cS \to \cS_0$ implies that $\cS_i$ is an analytic set  in $Q_i$ and $(\tau_i)|_{\cS} : \cS \to \cS_i$ is also finite. 
Furthermore, by construction $\cS_i$ meets $D_i$ along the union of analytic germs of $Z_i$
at the points from $\psi_i^{-1}(z_0)\cap \cS_i$. 
Assume by induction  that  for a germ of some holomorphic function $h_i$ on $Q_i$ the variety $\cS_i$ coincides with $h_i^*(0)$ and  the zero multiplicity of $h_i \circ \tau_i$ at general points of $E$
is $kl_im_i$ where $m_i$ is the degree of the finite morphism $\theta_i : V \to Z_i$ and $l_i=l-i$ (note that $l_i$ must be greater than zero since otherwise $Z_i$ is dense in $D_i$).

Hence the function $h_i/f^{km_i}$ has a holomorphic lift to the normal variety $X$. Furthermore, it is constant along each fiber of $\rho$ and it can be also pushed to a holomorphic 
function on the normal variety $Q$. This implies that $\delta_i$ factors through the Davis modification $\kappa : Q' \to Q_i$ of $Q_i$ along $D_i$ with ideal generated by $f^{km_i}$ and $h_i$.
By Lemma \ref{technical} for a general point $z'$ of $Z_i$
there is a Euclidean neighborhood with a local coordinate system $(u',v',\bar w')$ such that $D_i$ is given locally by $u'=0$, $Z_i$ by $u'=v'=0$, and 
up to an invertible factor
 $h_i$ can be presented locally as a product of $km_i$ factors of form $v'-u'd(u',\bar w')$. Hence by Lemma \ref{ref.l3} %Claim in Proposition \ref{ru.p2.4} 
 $\kappa$ is locally homogeneous at $z'$.
Thus $\delta_i$ is semi-finite and there is a simple modification $\sigma_{i+1} : Q_{i+1} \to Q_i$ as required in Lemma \ref{5.4}.

In order to finish induction it remains to prove the existence of a function $h_{i+1}$.  Note that over $z'$ the function $h_i\circ \sigma_{i+1}$ has zero multiplicity $km_i$ at general 
points of $D_{i+1}\subset Q_{i+1}$. Thus $h_{i+1}=h_i\circ \sigma_{i+1}/f^{km_i}$ is a regular function on the normal variety $Q_{i+1}$ that does not vanish at general points of $D_{i+1}$.
In particular $h_{i+1}$ vanishes only on the proper transform of $\cS_i$ which is, by construction,  $\cS_{i+1}$. This is the desired function and we are done.
\eproof 

%\brem\label{r.finite} Let $\cS$ and $\cS_i$ be as in the proof of Lemma \ref{5.5}. Note that we can require that $\cS$ contains a given point $x \in E$.
%Let $\cS'$ be the union of the components of $\cS$ passing through $x$. The induction in the proof of Lemma \ref{5.5} implies that when $Z_i$ is not dense in $Q_i$
%then the restriction of $\tau_{i+1}$ yields a finite holomorphic map $\cS \to \cS_{i+1}$. Thus, restricting this map further we get a finite holomorphic map
%$\cS' \to \cS_{i+1}'$ where $\cS_{i+1}'$ is  the union of irreducible components of $\cS_{i+1}$ passing through the point $\tau_{i+1} (x)$. That is,
%$\cS_{i+1}'$   a germ of an analytic hypersurface given by zeros of a holomorphic function in $Q_{i+1}$. Furthermore, we can require
%that $\cS'$ is contained in a principal hypersurface $W \subset X$ as in Remark \ref{ru.r5.5}.
%
%\erem

\bthm\label{5.6}
Let Notation \ref{5.1} and Convention \ref{5.2a} hold. Then $\C [X]^\Phi$ is finitely generated and $Q=X//\Phi$ coincides with the affine algebraic variety
$Q_n$ from Notation \ref{5.3} for some $n$.

\ethm

\bproof  By Proposition \ref{2.12} the sequence
\bdi Q_{n} & \rTo^{\sigma_{n}}& Q_{n-1}&\rTo^{\sigma_{i}}& \ldots &\rTo^{\sigma_2} & Q_1 &\rTo^{\sigma_1}& Q_0 \edi
of simple modifications from Notation \ref{5.3} cannot be extended indefinitely to the left. 
Hence by Lemmas \ref{5.4} and \ref{5.5} we can suppose that for some $n$ the image $\tau_n(E)$ is dense in $D_n$.
Suppose that $\rho : X \to Q$ is any partial quotient morphim over $Q_n$. Then $\delta_n(D)$ is dense in $D_n$ since $\tau_n$ factors through $\delta_n$.
Then by Lemma \ref{5.2} $Q$ is isomorphic to $Q_n$.  

Recall that $\C [X]^\Phi = \varinjlim \C [Q_\alpha]$ where $Q_\alpha$ is a partial quotient. Since every such a quotient over $Q_n$ is $Q_n$ itself
we see that $\C [X]^\Phi = \C [Q_n]$ which yields the desired conclusion.

\eproof

%Let $n=i+1$ in  Remark \ref{r.finite}. Then Theorem \ref{5.6} implies the following.
%
%\bcor\label{c.finite} Let the assumption of  Theorem \ref{5.6} hold. Then for every point $x \in E$ and $q=\rho (x) \in Q$ there exists 
%an analytic germ $\cS$ (resp. $\cS_0$) of an algebraic variety in $X$ through $x$ (resp. in $Q$ through $q$) of dimension $\dim Q -1$ such that
%$\cS_0$ is the zero locus of a holomorphic function and
%$\rho|_\cS : \cS \to \cS_0$ is a finite holomorphic map. Furthermore, if $x$ belongs to a principal hypersurface $W\subset X$ as  in Remark \ref{ru.r5.5}
%one can suppose that $\cS$ is contained in $W$. % and $\cS \setminus E=(\rho|_W)^{-1}(\cS_0\setminus D)$.
%
%\ecor

\bcor\label{5.7} Let $X$ be four-dimensional and let Notation \ref{5.1} hold without assuming
Convention \ref{5.2a}.  Suppose that $E$ admits a nonconstant morphism into a curve if and only if this
curve is a polynomial one. 

Then Convention \ref{5.2a} holds and thus  $\C [X]^\Phi$ is finitely generated. In particular Statement (1) of Theorem \ref{0.2} is true.

\ecor

\bproof By Remark \ref{5.2b} Convention \ref{5.2a} is true if $Z_0$ is not a point.
Assume it is a point $z_0$. Since $Q_0$ and $D_0=f^*(0)$ are smooth we can suppose that $z_0$ is locally a strict complete intersection
given by $f=g_1=g_2=0$. Note  that $g_i\circ \tau$ vanishes on $E$. Hence $(g_i/f)\circ \tau$ is a regular $G_a$-invariant
function on the normal variety $X$. This implies that $\psi : Q \to Q_0$ factors  through $\sigma_1 : Q_1 \to Q_0$ where $\sigma_1$ is the simple
modifications along $D_0$ with the ideal generated by $f,g_1$ and $g_2$. By construction $Q_1$ is smooth and the exceptional divisor $D_1 \simeq \C^2$. 
Thus we can replace the pair $(Q_0,D_0)$ with $(Q_1, D_1)$. If the center of $\psi$ after this replacement is still a point we observe
that it is again a strict complete intersection $f=h_1=h_2=0$ with $h_i=(g_i/f)\circ \sigma_1$. That is, the zero multiplicity of the lift $h_i$ (to $X$) on $E$ is less than 
the one of $g_i\circ \tau$.
Continue this procedure. As soon as one of these multiplicities becomes zero
the variety $Z_0$ becomes at least a curve. Thus we are done.

\eproof

Another consequence of the equality $Q=Q_n$ in Theorem \ref{5.6} is that the fiber of the morphism $Q \to B$ over $o$ is $D_n\simeq Z_{n-1}\times \C$. By the assumption
of Theorem \ref{0.2} $Z_{n-1}$ is a polynomial curve.
Hence we have the following.

\bcor\label{5.8}
Statement (2) of Theorem \ref{0.2} is true.
\ecor

\section{Proper actions}

\bdefi\label{6.0}  Recall that an action of an algebraic group $G$ on a variety $X$ is proper
if the morphism $G\times X \to X \times X$ that sends $(g,x) \in G \times X$ to $(x, g.x) \in X \times X$ is proper. 
\edefi

It is well-known that every proper $G_a$-action is automatically free. Some geometrical properties of such actions are described 
below.

%\bnota\label{6.0a} In this section $X$ is a normal affine algebraic variety equipped with a proper $G_a$-action $\Phi$
%and $\rho : X \to Q$ is the quotient morphism in the category of affine algebraic varieties.
%\enota

\bprop\label{6.1} Let  $X$ be a normal affine algebraic variety equipped with a proper $G_a$-action $\Phi$  and let $C_i$ be a sequence of general orbits
of $\Phi$. Suppose that $C\subset X$ is a closed curve such that for every $c \in C$ and every neighborhood $V\subset X$ of $c$ (in the standard topology) there exists
$i_0$ for which $C_i$ meets $V$ whenever $i \geq i_0$. Then $C$ is an orbit of $\Phi$,
i.e. $C \simeq \C$ (in particular it is smooth and connected).
\eprop

\bproof 
Assume that $C$ meets two disjoint orbits $C'$ and $C''$ of $\Phi$ and choose points
$c_i', c_i'' \in C_i$ such that $c_i' \to c' \in C'\cap C$ and $c_i''\to c''\in C''\cap C$ as $i \to \infty$. Note that $c_i''=t_i.c_i'$. This sequence of numbers $\{t_i \}$ in the group $\C_+$
goes to $\infty$. Indeed, otherwise switching to a subsequence we can suppose that $\lim_{i \to \infty } t_i =t$ and by continuity
$c''=t.c'$ which is impossible since the last two points are in distinct orbits. But this implies that the preimage (in $G\times X$) of a small neighborhood
of the point $(c',c'') \in X\times X$ contains points of form $(t_i,c_i')\in \C_+\times X$ going to infinity contrary to properness.

Hence $C$ meets only one orbit $O$ of $\Phi$. 
Hence  $C\setminus O = \emptyset$ since otherwise $C$ meets other orbits.
Therefore, $C= O\cap C$ which implies that $C=O\simeq \C$
% Being a limit of $\Phi$-invariant curves $C$ is also $\Phi$-invariant.
%Hence $O$ cannot be a point since otherwise the action is trivial on $C$ which therefore contains infinite number of orbits
\footnote{This implies, in particular, that $\Phi$ is free.} and we are done. 
\eproof

\bcor\label{6.2} Suppose that for some $q \in Q$ the fiber $\rho^{-1}(q)$ is a curve. Then $\rho^{-1}(q)\simeq \C$.
\ecor

Let us fix notation for the rest of this section.

\bnota\label{ref.n6} We suppose that
$B$ is a unibranch germ of a smooth complex algebraic curve at point $o=f^*(0)$ (where $f \in \C [B]$)
and  $\varphi : X \to B$ is a morphism from a complex factorial affine algebraic variety $X$ equipped with a $G_a$-action $\Phi$
which preserves each fiber of $\varphi$ and such that the fiber $E=\varphi^*(o)$ is reduced.

We suppose also that Convention \ref{5.2a} is true for some partial quotient morphism
$\rho : X\to Q$  of the action $\Phi$ and a morphism $\psi : Q\to Q_0$
into a smooth affine algebraic variety $Q_0$ over $B$ described in Notation \ref{5.1}.
Recall that in this case by Theorem \ref{5.6} $\C [X]^\Phi$ is finitely generated and $Q=X//\Phi$ coincides with an affine algebraic variety
$Q_n$ from Notation \ref{5.3} for some $n$, and, particular, in those notations $\tau_n : X \to Q_n$ is the quotient morphism.
Furthermore, let $D_n$ be the divisor in $Q_n$ over $o \in B$. Then it follows from the construction of $Q_n$ that $D_n \simeq Z_{n-1} \times \C$
and from the proof of Theorem \ref{5.6} that $\tau_n (E)$ is dense in $D_n$.

\enota

\bprop\label{6.3} Let  %the assumption of Theorem \ref{5.6} and Notation \ref{5.3} hold (i.e. $Q=Q_n\supset D_n$ and $\tau_n(E)$ is dense in $D_n=Z_{n-1}\times \C$) in addition to Notation \ref{6.0a}
$\alpha : E \to R$ be the quotient morphism of the restriction of $\Phi$ to $E$. Then there is a natural affine
modification $\beta : R\to D_n$ over $Z_{n-1}$.

\eprop

\bproof 
The density of $\tau_n (E)$ in $D_n$ implies that for a general point $q \in D_n$ its preimage $\tau_n^{-1}(q)$ is a curve.
By Corollary \ref{6.2} this curve is isomorphic to $\C$.

By the universal property of quotient morphisms the map $\tau_n|_E : E \to D_n$ factors through $\alpha$, i.e. there exists 
$\beta : R \to D_n$ for which $\tau_n|_E= \beta \circ \alpha$.  Since $\tau_{n-1}$ (which maps $E$ onto $Z_{n-1}$)  factors through $\tau_n$
we see that $\beta$ is a morphism over $Z_{n-1}$. Note that $\beta^{-1}(q)$ is a point
since $\tau_n^{-1}(q)$ is a connected curve. Thus $\beta$ is birational.
Any birational morphism of affine algebraic varieties is an affine modification
 \cite{KaZa} and we are done.

\eproof

\brem\label{6.5}  Proposition \ref{6.3} implies that $\beta$ is an isomorphism over a Zariski dense open subset $D_n^*$ of $D_n$.
Furthermore removing a proper subvariety from $D_n^*$ one can suppose that for every point in $\beta^{-1}(D_n^*)$ the fiber
of $\alpha$ over this point is an orbit of $\Phi$.
Let $Q'=(Q\setminus D_n)\cup D_n^*$. Note that by construction every fiber of $\rho$ over point of $Q_*$ is an orbit of $\Phi$
and the codimension of $Q \setminus Q'=D_n \setminus D_n^*$ in $Q$ is at least 2.

\erem

Recall that when $X$ is four-dimensional and there is no nonconstant morphism from $E$ into
any non-polynomial curve then by Remark \ref{5.2b} Convention \ref{5.2a} holds,  $Z_{n-1}$ is a polynomial curve, and thus the normalization of $D_n=Z_{n-1}\times \C$ 
is $\C^2$.  In particular, if $E$ is normal then  $R$ is normal (being the quotient
space of normal $E$) and  we can mention the following interesting fact (which won't be used later).

\bcor\label{6.4} Let the assumption of Corollary \ref{5.7} hold and $E$ be normal. Then  $R$ is an affine modification of $\C^2$.
%Under the assumption of Corollary \ref{5.7} $R$ is an affine modification of $\C^2$.
Furthermore, for a factorial $E$ one has $R \simeq \C^2$ and normalization of morphism $\beta$ is a birational morphism $\beta' : \C^2 \to \C^2$.  
\ecor

\bproof Since $\beta$ is a morphism over the curve $Z_{n-1}$ we see that $\beta' : R \to \C^2$ is a morphism
over the normalization $Z_{n-1}^{\rm norm}\simeq \C$ of $Z_{n-1}$.  Note that by Corollary \ref{fv.c1} $R$ is factorial being the quotient space
of the factorial threefold $E$.  This implies that all fibers of the $\C$-fibration $R \to Z_{n-1}^{\rm norm}$ are reduced and irreducible.
Hence $R \simeq \C^2$ and we are done.

\eproof

\bprop\label{6.6} 
Let $\beta$ be as in Proposition \ref{6.3} and suppose that one of the following conditions is satisfied:

{\rm (a)} $\beta$ is finite, $E$ is smooth outside codimension 2, and every fiber of $\alpha : E \to R$ is a curve
with a possible exception of a finite number of such fibers;

{\rm (b)} $\beta$ is quasi-finite and  the assumption of Theorem \ref{0.2} (3) holds (in particular, $X$ is Cohen-Macaulay, $\Phi$ is proper, each fiber $X_b$ of $\varphi : X \to B$ (including $E$)
is normal, and, being the image of $E$, the curve $Z_{n-1}$ is a polynomial one).

Then $Z_{n-1}$ (and therefore $D_n =Z_{n-1}\times \C $) is smooth and $\beta$ is an isomorphism. 

%be finite and $E$ be smooth outside codimension 2. Suppose also that every fiber of $\alpha : E \to R$ is a curve
%with a possible exception of a finite number of such fibers. 
%Then $Z_{n-1}$ (and therefore $D_n =Z_{n-1}\times \C $) is smooth and $\beta$ is an isomorphism. 

%Furthermore, suppose that the assumption of Theorem \ref{0.2} (3) holds (in particular, $Z_{n-1}$ is a polynomial curve and $E$ is normal). Then the statement holds even in the case when
%$\beta$ is quasi-finite.
\eprop

\bproof

The assumptions on $\alpha$ and $\beta$ imply that there is a finite subset $M\subset D_{n}$ such that for every $q \in D_n\setminus M$ the preimage
$\tau_n^{-1}(q)$ is a curve. 

Recall that $X$ is Cohen-Macaulay which implies that the morphism $\rho$ is faithfully flat  over $Q \setminus M$ (e.g., see \cite[Chap. 2 (3.J) and Th. 3]{Mat} and \cite[Th. 18.6]{Ei}).
Hence $\Phi$ is locally a translation over $Q\setminus M$  by  \cite[Theorem 2.8]{DeFiGe} (see also \cite[Theorem 1.2]{DeFi}).
%and $Z_{n-1} \simeq \C$
This implies that $Z_{n-1}$ is smooth since otherwise $D_n\setminus M$ and 
therefore $(D_n\setminus M) \times \C\subset E$ are not smooth outside codimension 2. 
It remains to note that, being birational and finite, $\beta$ is an isomorphism by the Zariski Main theorem.

For the second statement consider the normalization  $\kappa : R^{\rm norm} \to D_{n}^{\rm norm}\simeq \C^2$ of $\beta$. Since $\beta$ is birational by Proposition \ref{6.3} and quasi-finite
by the assumption, the morphism $\kappa$ is an embedding. Since both   $R^{\rm norm}$ and  $D_{n}^{\rm norm}$ are affine we see that
 $D_{n}^{\rm norm}\setminus \kappa (R^{\rm norm})$ is either empty or a Cartier divisor. Since $ D_{n}^{\rm norm}\simeq \C^2$ such a divisor must be given by zeros of a regular function
 $g$. Then $g$ admits a regular lift to the normal variety $E$ that does not vanish. This contradicts the assumption that $E$ does not admit nonconstant morphisms  into a non-polynomial curve.
Thus $\kappa (R^{\rm norm} )= D_{n}^{\rm norm}$ and as before the assumption on $\alpha$ implies  
 that  there is a finite subset $M\subset D_{n}$ such that for every $q \in D_n\setminus M$ the preimage
$\tau_n^{-1}(q)$ is a curve. Hence we can repeat the previous argument and get the desired conclusion.

\eproof

\brem\label{ru.re7.1} (1) Let $Q'$ be from the Remark \ref{6.5}.  Note that exactly the same argument as in Lemma \ref{6.6} for the set $Q\setminus M$
implies that the action $\Phi$ over $Q'$ is locally trivial\footnote{That is, $\rho^{-1}(Q')$ can be covered
by  $\Phi$-stable affine open subsets on each of which the action admits an equivariant trivialization.} even without the assumption that $\beta$ is quasi-finite. 

(2) Actually, following \cite{Ka02} one can extract the Sathaye's theorem \cite{Sa} from the argument in the proof of Proposition \ref{6.6}. 
Indeed, under the assumption of  the Sathaye's theorem  $D_{n}=\C^2$ and thus  $Z_{n-1} \simeq \C$.
Since $Z_{n-1}$ is smooth every polynomial curve $Z_i$ is smooth by Lemma \ref{3.3}
and therefore  $Z_{i} \simeq \C$. Hence every $D_{i}= Z_{i-1} \times \C$ is the plane and by the Abhyankar-Moh-Suzuki theorem
$Z_{i}$ can be viewed as a coordinate line in $D_{i}$. It is a straightforward fact that a modification $\sigma_{i+1} : Q_{i+1} \to Q_{i}$ with
such center $Z_i$ and divisor $D_i$ leads to $Q_{i+1}$ isomorphic to $Q_i$ which implies that $Q_n \simeq B\times \C^2$.

(3) The assumption on properness of $\Phi$ is crucial in Proposition \ref{6.6}. In the absence of properness Winkelman \cite{Wi} constructed a free
action on $X=\C_{x_1,x_2,x_3,x_4}^4$ such that the quotient $Q$ is isomorphic to the hypersurface in $\C^4_{x_1,u,v,w}$ given by
$$x_1w=v^2-u^2-u^3.$$  Both $X$ and $Q$ are considered over the curve $B=\C_{x_1}$ and in particular the zero fiber $D$ of the morphism $Q \to B$
is given by $v^2-u^2-u^3=0$ in $\C^3_{u,v,w}$. That is, $D$ is not a plane but it is the product of $\C$ and
$Z_{n-1}$ where $Z_{n-1}$ is a polynomial curve with one node as singularity (in accordance with Theorem \ref{0.2} (2)).
\erem

\section{ Some facts about fibered products}

In order to establish quasi-finiteness required in Proposition \ref{6.6} we need some auxiliary facts.

\bnota\label{fp.n0}    Let $\rho : X \to Q$ be a dominant morphism of normal quasiprojective algebraic varieties and $S$ be an irreducible closed subvariety of $X$
such that for a Zariski dense open subset $S'$ of $S$ the restriction $\rho |_{S'}: S' \to Q$ is quasi-finite. Suppose that $k = \dim X - \dim S$ is the
dimension of general fibers of $\rho$.
\enota

%The next fact follows from the semi-continuity theorem \cite[Chapter III, Theorem 12.8]{Har}  \cite[Chapter I, Section 6.3, Corollary]{Sha} .

\blem\label{fp.l-1}
Consider an irreducible germ $\Gamma \subset S$ of a general curve through $s \in S$ (in particular, $\Gamma \setminus s \subset S'$ and
the dimension of $\rho^{-1}(\rho (\Gamma  \setminus s))$ is $k+1$). Then
the dimension of the fiber $F$ in the closure of  $\rho^{-1} (\rho (\Gamma \setminus s ))$ (in $X$) over $s$ is $k$.
\elem

\bproof  By the  semi-continuity theorem  \cite[Chapter I, Section 6.3, Corollary]{Sha} the dimension of $F$ is at least $k$.
On the other hand it cannot be greater than $k$ since otherwise it is not contained in the closure of the $(k+1)$-dimensional variety $\rho^{-1}(\rho (\Gamma  \setminus s))$.

\eproof

Note that for a fixed $s$ this fiber $F$ may in general depend on the choice of $\Gamma$, i.e. $F=F(\Gamma )$. 

\bdefi\label{fp.d1} Suppose that $F_0(\Gamma )$ is the component of $F(\Gamma )$ containing $s$. We say that the data in Notation \ref{fp.n0} satisfies Condition (A) if

\noindent {\em for every $s \in S$ the set $\bigcup_\Gamma F_0(\Gamma)$
consists of a finite number of $k$-dimensional varieties.}

\edefi

\bexa\label{fp.ex1} (1) If $X=S$ then Condition (A) holds automatically.

(2) Suppose that $\rho : X \to Q$ is the quotient morphism of a proper $G_a$-action $\Phi$ and $S$ is a hypersurface of $X$ such that $\rho (S)$ is  Zariski dense in $Q$.
Then for a general germ $\Gamma$ as in Definition \ref{fp.d1} and every point $s' \in \Gamma\setminus s$ the fiber $\rho^{-1}(\rho(s'))$ is a union of orbits of $\Phi$.
The set  $\bigcup_\Gamma F_0(\Gamma)$ contains, of course, the orbit of $\Phi$ through $s$ and by Proposition \ref{6.1} it contains nothing else. Thus Condition (A) holds.
\eexa

\bprop\label{fp.p1}    Let Notation \ref{fp.n0} hold and Condition (A) be satisfied.
Suppose that $\tX$ is the irreducible component of $X\times_QS$ such that it contains $\tS =\{ (s,s)| s\in S \}$.
%the restriction of  the natural projection $X\times_QS \to X$  to any of these components is dominant. 
Denote  by  $\lambda : \tX \to S$ the restriction  to $\tX$ of the natural projection.

Then for every irreducible subvariety $P \subset S \setminus S'$ the dimension of $\lambda^{-1} (P)$  coincides with $m +k$ where $m=\dim P$.

\eprop

\bproof 

Suppose that for $l \geq 0$ the subvariety $P(l)\subset P$ consists of points $s \in P$ for which $\dim \lambda^{-1}(s)=k+l$. Our aim is to show that
$\dim P(l) \leq m-l$ and thus  $\dim \lambda^{-1}(P) = m+k$.

Let $s \in S, \Gamma ,  F$ be as in Definition \ref{fp.d1} and let
the closure of $\rho^{-1} (\rho (\Gamma \setminus s))\cap S$ 
meets $F$ at  points
$s_1, \ldots , s_r$ (i.e. each irreducible component of $F$ passes at least through one of these points).
Note that for a fixed $s$ this set  $\{ s_1, \ldots , s_r \}$ may in general depend on the choice of $\Gamma$, i.e. $s_i=s_i (\Gamma )$. 
Suppose that $M(s)$ is the union $\bigcup_i\bigcup_\Gamma s_i(\Gamma)$.  By definition of $\tX$ for every $s' \in \Gamma \setminus s$ one has $\lambda^{-1}(s')\simeq \rho^{-1}((\rho (s'))$. Hence
by Lemma \ref{fp.l-1} and Condition (A),   
the variety $\lambda^{-1}(s) =\bigcup_\Gamma F(\Gamma )$ has dimension at most $\dim M(s)+k$.
In particular $P(l)$ is contained in the set $$\{ s \in P | \dim M(s)\geq l \}$$ and we need to estimate the dimension of this last set.

For this we need to use the Hironaka flattening theorem \cite{Hi} which implies the existence  of a proper birational morphism $\hQ \to Q$ such that
for the union $\hS$ of irreducible components of $S\times_Q \hQ$ with the dominant projections to $S$ the induced morphism $\hS \to \hQ$ is quasi-finite.
We have the following commutative diagram 
\bdi
\hS &\rTo^{\beta}&S\\
\dTo>{\theta}&&\dTo>{\gamma}\\
\hQ&\rTo^{\alpha}&Q\\
\edi
where the morphisms $\alpha$ and $\beta$ are birational and the morphism $\theta$ and $\gamma |_{S'}$ are quasi-finite.

Let $P=\bigcup_i P_i$ where $P_i$ are disjoint (not necessarily closed) irreducible varieties such that for any irreducible component $V$ of $\beta^{-1}(P_i)$ 
all fibers of the morphism
$\beta|_{V} : V \to P_i$ are of the same dimension. 
Varying $i$ and $V$ consider the collection $C_1, C_2, \ldots $ of irreducible subvarieties of $\hQ$ that are of the form $C_j=\overline{\theta (V)}$.  Observe
that stratifying the varieties $P_i$'s further one can suppose that  for indices $j \ne l$ the varieties $C_j$ and $C_l$ are either equal or disjoint. 
Denote by $I$ the set of all pairs $(i,j)$ for which there exist $V\subset \beta^{-1}(P_i)$ with $\theta (V)=C_j$
and put $C_{ij}=V$  (in particular $\dim C_j =\dim C_{ij}$
because of quasi-finiteness of $\theta$).

Let  $n_{ij}$ be the dimension of a fiber $F_{ij}\subset \beta^{-1}(s)$ of the morphism $\beta |_{C_{ij} } : C_{ij} \to P_i$. 
In particular
$\dim P_i = \dim C_j - n_{ij} \leq m$. Hence for $n_j=\min_i \{ n_{ij}| (i,j) \in I \}$ one has $\dim P_i \leq m+n_j-n_{ij}$ for any $(i,j) \in I$.

Let $s \in P_i$, i.e. $\beta^{-1} (s) =\bigcup_j F_{ij}$. For every  where $(tj) \in I$  we put $F^{t}_{ij}=\theta^{-1} (\theta (F_{ij}))\cap C_{tj}$, i.e. $\dim F^{t}_{ij}=n_{ij}$ because of quasi-finiteness.
Note that $$\dim \beta (F^{t}_{ij})=\max (0, n_{ij}-n_{tj}) \leq n_{ij} -n_j\leq \max_j (n_{ij} -n_j).$$ Thus $l:=\dim \beta (\bigcup_{t}F^{t}_{ij}) \leq \max_j (n_{ij} -n_j)$ while $\dim P_i \leq m- \max_j (n_{ij} -n_j)\leq m-l$. 

Since we started with a point $s \in P_i$ for the desired inequality $\dim P(l) \leq m-l$ it suffices to show that $M(s) \subset \bigcup_j\beta (\bigcup_{t}F^{t}_{ij})$
(since the last set is contained in $P_i$).
But this follows from the fact that the proper transform of $\Gamma$ (mentioned in the definition of $M(s)$) in $\hS$ meets $\beta^{-1}(s)$ and therefore
meets some $F_{ij}$. By quasi-finiteness 
the proper transform (in $\hS$) of every component in the closure of $\rho^{-1} (\rho (\Gamma \setminus s))\cap S$ meets some $F^{t}_{ij}$. This implies that this closure passes through a point in 
$ \beta (F^{t}_{ij})$ which yields $M(s) \subset \bigcup_j\beta (\bigcup_{t}F^{t}_{ij})$ and we get  $\dim \lambda^{-1} (P) \leq m+k$, i.e.  $\dim \lambda^{-1} (P)\leq \dim \rho^{-1}(\rho (P))$.

\eproof

\bcor\label{fp.c1} If Condition (A) holds and $\rho^{-1}(\rho (P))$ is at least of codimension  2 in $X$ then  $\lambda^{-1} (P)$
has codimension at least  2 in $\tX$.
\ecor

%\bcor\label{fp.c2} Let $\tS$ be the union of components of $S\times_QS$ whose natural projections to $S$ are dominant (i.e. $\tS \subset \tX$).
%Then  morphism $\lambda|_\tS : \tS \to S$ is quasi-finite.

%\ecor

%\bproof Apply Proposition \ref{fp.p1} to the case of $X =S$ and $P$ is a singleton.

%\eproof

\section{Main Theorems}

\bnota\label{7.1} 
Up to Section \ref{10.5a} below we adhere to the assumptions of Proposition \ref{6.3} and, in particular, Notation \ref{5.1}. That is, $B$ is the germ
of a smooth algebraic curve at point $o=f^*(0)$, $\varphi : X \to B$ is a morphism of factorial varieties, $\Phi$ is a $G_a$-action on $X$,
$Q=Q_n$ is the categorical quotient $X//\Phi$ in the category of affine algebraic varieties, 
$\tau_n=\rho$, and $\rho (E)$ is dense in $D_n=:D$.  Since $\beta$ is birational by Proposition \ref{6.3}, there is a Zariski dense open $R^* \subset R$
such that the restriction of $\beta|_{R^*}: R^* \to D\subset Q$ to $R^*$ is an embedding. 
%Furthermore since $\beta|_R$ is a morphism over a polynomial curve $Z_{n-1}$
%we can suppose that $R^*$ is isomorphic to a product $Z^* \times \C$ where $Z^*$ is a smooth rational curve.

Suppose also (as in Theorem \ref{0.2} (3)) that the restriction of 
our proper $G_a$-action $\Phi$ to $E$ is a translation. In particular,
the quotient morphism $\alpha : E \to R$ admits a section whose image in $E$, by abuse of notation, will be denoted by $R$.
The image of $R^*$ in this section will keep notation $R^*$ as well
%Let $\hR^*$ the preimage of $R^*$ in $\hR$
(i.e. $\rho |_{R^*} : R^* \to D$ is an embedding).

The section $R$ coincides with the zero locus of a regular function on $E$ which has simple zeros at every point of $R$.
Extend this function to a regular function $h$ on $X$ and consider the zero locus  $h^{-1}(0)=:S  \subset X$ of this extension.
\enota

\blem\label{l8.7} The function $h$ from Notation \ref{7.1} can be perturbed so that for $S^*=S\setminus R$ and $Q^*=Q\setminus D$ the restriction
$\rho|_{S^*}: S^* \to Q^*$ is finite. In particular, for $S' =S^*\cup R^*$ the restriction
$\rho|_{S'}: S' \to Q$ is quasi-finite.

\elem

\bproof Recall that $X\setminus E \simeq Q^*  \times \C$. Let $w$ be a coordinate on the second factor, i.e.
$h|_{X\setminus E}$ can be treated as a polynomial in $w$ with coefficients in $\C [Q^*]$. Let $n$ be the degree of this polynomial
and $f$ be as in Notation \ref{5.1}. Then for sufficiently large $k$ the function $h+f^kw^{n+1}$ is regular on $X$
and the restriction of $h+f^kw^{n+1}$ to $E$ has only simple zeros on its zero locus $R$. This $h+f^kw^{n+1}$ yields the desired perturbation.

\eproof

\bprop\label{mt.p1}
Suppose that $\tX$ is the irreducible component of $X\times_QS$ such that it contains $\tS=\{ (s,s) | s \in S\}$. %the restriction of  the natural projection $X\times_QS \to X$  to this component is dominant. 
Then $\tX$ is naturally isomorphic to $S \times \C$ with $\tS$ being the section of the natural projection $\lambda : \tX \to S$ .
\eprop

\bproof Let $S'$ be as in Lemma \ref{l8.7} and $s \in S'$.  Note that $\rho^{-1}(\rho(s))$ is a disjoint union of orbits of $\Phi$ because of quasi-finiteness of $\rho|_{S'}$.
One of these orbits $O_s$ passes through $s$.  By construction  $\lambda^{-1}(\lambda (s)) \subset \rho^{-1}(\rho(s))\times s$.
Since $O_s \times s$ is the only component of $\rho^{-1}(\rho(s))\times s$ that meets $\tS$ we see that  $\lambda^{-1}(\lambda (s))=O_s \times s$
which is an orbit of the free action $\tilde \Phi$ on $\tX$ induced by $\Phi$ (in particular $\tilde \Phi$ preserves the fibers of $\lambda$). 
% Note that the proper action $\Phi$ on $X$ induces a free action $\tilde \Phi$ on $\tX$ such that $\tilde \Phi$ preserves
%the fibers of the natural projection $\lambda : \tX \to S$.
%Let $S'$ be as in Lemma \ref{l8.7}. Then because of quasi-finiteness of $\rho|_{S'}$ every fiber $\lambda^{-1}(\lambda (s)), \, s \in S'$
%coincides with $\rho^{-1}(\rho(s))\times s \simeq \C \times s$ which is an orbit of $\tilde \Phi$ (since $\rho^{-1}(\rho(s))$ is an orbit of $\Phi$).
The variety $\tS'=\{ (s,s) | s \in S'\}$ is a section of  $\tilde \Phi$ over $S'$. % and, furthermore, it meets every fiber $\lambda^{-1}(\lambda (s)), \, s \in S'$ once.
Hence $\lambda^{-1}(S')$ is naturally isomorphic to $S' \times \C$ because the action $\tilde\Phi$ is free.
By construction  $S \setminus S'$ is of codimension at least 2 in $S$. Hence the same is true for the subvariety  $(S \setminus S')\times \C$ in $S\times \C$. 
By Corollary \ref{fp.c1} and Example \ref{fp.ex1} (2) $\tX \setminus \lambda^{-1}(S')$ is of codimension at least 2 in $\tX$. By the Hartog' theorem
the isomorphism  $\lambda^{-1}(S') \simeq  S' \times \C$ extends to an isomorphism $\tX \simeq S \times \C$ since both varieties are affine.
This is the desired conclusion.
\eproof

\blem\label{mt.l1} For $\tX$  from Proposition \ref{mt.p1} the restriction $\kappa : \tX \to X$ of the natural  projection $X\times_QS \to X$ is quasi-finite.

\elem

\bproof The restriction of $\kappa$ to $\lambda^{-1}(S')$ is quasi-finite because of quasi-finiteness of $\rho|_{S'}$.
Since $\tX \setminus \lambda^{-1}(S')=(S\setminus S') \times \C=R \times \C$ it suffices to check that the restriction of $\kappa$ yields a quasi-finite  morphism $R \times \C \to E$.
The last map is an isomorphism since $R$ is a section of $\Phi|_E$ and we are done.

\eproof

\bnota\label{10.1} By the Grothendieck version of the Zariski Main theorem (e.g., see \cite[Th\'eor\`eme 8.12.6]{EGA}) 
there is an embedding $ \tX \to \hX$ of $\tX$ into an algebraic variety $\hX$ such that $\kappa$ can be extended to
a finite morphism $\chi : \hX \to X$. Note that the finiteness implies that $\hX$ is affine since $X$ is. Furthermore, replacing if necessary
$\hX , \tX$, and $S$ with their normalizations we suppose that these varieties are normal.
 
Denote the lift of $f \in \C [B]$ from Notation \ref{7.1} to $X$ (resp. $\tX$, resp. $\hX$) by
the same letter $f$ (resp. $\tilde f$, resp. $\hat f$). Let $\nu$ be the locally nilpotent vector field on $X$ associated with 
the action $\Phi$. Since $f \in \Ker \nu$ the field $f^k \nu$ is also locally nilpotent for every $k>0$ and it is associated
with a $G_a$-action $\Phi_k$ on $X$ which has the same quotient morphism $\rho : X \to Q$. The lift of $\nu$ to
$\tX$ (resp. $\hX$) will be denoted by $\tilde \nu$ (resp. $\hat \nu$). By construction $\tilde \nu$ (and, hence, $\tilde f^k \tilde \nu$) is locally nilpotent on $\tX$
and $\hat \nu$ is a rational vector field on $\hX$. Since $\hX \setminus \hat f^{-1}(0)\subset \tX$ we see that for sufficiently
large $k$ the field $\hat f^k\hat \nu$ is locally nilpotent on $\hX$. We denote by $\hat \Phi_k$ the $G_a$-action associated with the last field.

\enota

\bprop\label{10.2} The action $\hat \Phi_k$ admits a quotient morphism $\hat \rho : \hX \to \hQ$ in the category of affine algebraic varieties
such that the following commutative diagram holds
\bdi
\hX &\rTo^{\hat \rho}&\hQ\\
\dTo>{\chi}&&\dTo>{\tau}\\
X&\rTo^{\rho}&Q\\
\edi
where $\tau$ is a finite morphism.

\eprop

\bproof

For every function $g \in \C [\tX]^{\tilde \Phi}$ and sufficiently large $m>0$ the function $\tilde f^m g$ has a regular extension
to $\hX$, i.e. it can be viewed as a function from $\C [\hX]^{\hat \Phi_k}$. Since by construction $\hX \setminus \hat f^{-1}(0)= \tX \setminus \tilde f^{-1}(0) = \lambda^{-1} ( S\setminus R)$
where $\lambda : \tX \to S$ is the quotient morphism of $\tilde \Phi$ (by Proposition \ref{mt.p1})
we see now that a finite number of functions from $\C [\hX]^{\hat \Phi_k}$ separates the orbits of $\hat \Phi_k$ contained in $\lambda^{-1}(S \setminus R)$. 
Consider a partial quotient morphism $\hat \rho : \hX \to \hQ$ whose coordinates include these functions. Furthermore, since the lift of every function from $\C [Q]$
to $\hX$ is a function from $\C [\hX]^{\hat \Phi_k}$ we can include the lifts of the coordinate functions of $\rho$ into the set of coordinate functions of $\hat \rho$.
Then the commutative diagram as before make sense.

% and since  $\rho (E) \subset \tau \circ \hat \rho (\hE )$, where $\hE =\hat f^{-1}(0)$,
%we see that  $\tau \circ \hat \rho (\hE )$ is dense in $D$. 

Let us show that $\tau$ is finite. Indeed,
by finiteness of $\chi$ every function $g$ from $\C [\hX]^{\hat \Phi_k}$ is integral over $\C[X]$.
That is, $g$ is a root of an irreducible monic polynomial $P(g):=g^n+a_{n-1}g^{n-1}+ \ldots +a_1g +a_0$ where each $a_i \in \C [X]$.
Since $g$ is $\hat \Phi_k$-invariant it satisfies also an equation $g^n+a_{n-1}^tg^{n-1}+ \ldots +a_1^tg +a_0^t=0$ where $a_i^t$
is the result of the action of $t \in \C_+$ (induced by $\hat \Phi_k$) on $a_i$. However, the minimal polynomial $P$ should be unique which implies
that $a_i$ are $\hat \Phi_k$-invariant and thus $\C [\hX]^{\hat \Phi_k}$ is integral over $\C [Q]$. That is, every function from $\C [\hQ ]$ is integral over $\C [Q]$
and hence $\tau$ is finite.

This leads to the commutative diagram
\bdi
\hE &\rTo^{}&\hD\\
\dTo>{}&&\dTo>{}\\
E&\rTo^{}&D\\
\edi
where $\hD :=\tau^{-1}(D)$,  $\hE =\hat f^{-1}(0)$, and the vertical morphism are finite. 
Since for a general point $q \in D$ its preimage $\rho^{-1} (q)$ is a curve in $E$ the diagram implies that for a general point $\hat q \in \hD$ its preimage $\hat \rho^{-1} (\hat q )$
is  a curve in $\hE$.
Thus, by Theorem \ref{1.4} $\hat \rho$ is a quotient morphism which is the desired conclusion.

\eproof

\blem\label{10.3} Let $Q'=Q^* \cup D^*$ be as in Remark \ref{ru.re7.1}. % (i.e. $S'=S\cap \rho^{-1}(Q')$ where $S'$ is as in Lemma \ref{l8.7}).
Then 

{\rm (i)} $\hX \setminus (\rho \circ \chi)^{-1} (Q')$ has codimension at least 2 in $\hX$;

{\rm (ii)} for every point in $Q'$ there is a Zariski neighborhood $V \subset Q'$ such that $ (\rho \circ \chi)^{-1} (V)$ is naturally isomorphic to $\hat V \times \C$ where $\hat V =\tau^{-1}(V)$. 
\elem

\bproof Recall that $\rho^{-1}(Q\setminus Q') =E\setminus (\rho|_E)^{-1}(D^*)$ is of codimension at least 2 in $X$. Hence (i) is a consequence of finiteness of $\chi$.

By Remark \ref{ru.re7.1} $\Phi$ is a locally trivial over $Q'$. In particular, for every point in $Q'$ there is a Zariski neighborhood $V \subset Q'$ such that
$\rho^{-1}(V)$ is naturally isomorphic to $V \times \C$ where the projection $\rho^{-1}(V)\to V$ is the restriction of $\rho$.
Note that the restriction of the commutative diagram in Proposition \ref{10.1} yields another commutative diagram
\bdi
\hW &\rTo^{}&\hV\\
\dTo>{}&&\dTo>{}\\
\rho^{-1}(V) &\rTo^{}&V\\
\edi
where $\hW= (\rho \circ \chi)^{-1} (V)$. Hence we have $\hW\simeq \hV \times \C$ since $\rho^{-1}(V) \simeq V \times \C$. This yields (ii).

\eproof

\bprop\label{10.4} Let  $\hS$ be the closure of $\tS$ in $\hX$. Then $\hS$ is a section of the quotient morphism $\hat \rho : \hX \to \hQ$.

\eprop

\bproof If suffices to construct a $\hat \Phi_k$-equivariant morphism $\psi : \hX \to \hS$ from $\hX$ to $\hS$. Indeed, then by the universal properties of quotient morphisms
this morphism factors through $\hQ$ and, therefore, $\hat \rho|_\hS : \hS \to \hQ$ is invertible. 

Choose as the restriction of such a $\psi$ to $\tX$ the natural projection $\lambda : \tX=S \times \C \to S$. Note that for $\hV\simeq V \times \C$ from Lemma \ref{10.3}(ii)
$\lambda$ agrees on $\tX \cap \hV$ with the natural projection $\hV \to V$. Hence
we can extend $\psi$ to $\tX \cup \chi^{-1}( Q')$. 
Now because of Lemma \ref{10.3} (i) and the Hartogs' theorem this restriction extends to $\hX$ and we are done.

\eproof

\brem\label{future} The fact that the morphism $\alpha : E \to R$ has a section may be made weaker for Proposition \ref{10.4}.
It suffices, say, to require that there exists is a reduced irreducible principal effective divisor $T=h^*(0)$ in $E$ such that  the restriction $\alpha|_T : T \to R$ is surjective and quasi-finite.
Then one can extend $h$ to a regular function on $X$ and the same proofs as before imply that the result remains valid.

\erem

\bcor\label{10.5} The morphism $\beta : R \to D$ is quasi-finite. 

\ecor

\bproof  Let $\tR$ be the preimage of $R \subset S$ in $\tS$.  Note that $\tR$ is a subvariety of $\tX \subset \hX$. On the other hand since $\hQ$ is isomorphic to 
$\hS \supset \tS$ by Proposition \ref{10.4}, we can treat $\tR$ as a subvariety of $\hQ$. Furthermore, by construction the quotient morphism
$\hat \rho$ yields an automorphism of $\tR$. Hence the restriction of the morphism $\tau \circ \hat \rho = \rho\circ \chi$ to $\tR\subset \hX$ yields a quasi-finite morphism $\tR \to D$ (where
$\tau, \chi$, $\hat \rho$ are from Proposition \ref{10.2}).
Since morphism $ \rho\circ \chi |_\tR$ factors through $\beta =\rho|_R : R \to D$ the latter is also quasi-finite and  we are done.
\eproof

 \subsection{Proof of Theorem \ref{0.2}}\label{10.5a}   Note that Convention \ref{5.2a} is valid by Remark \ref{5.2b}. Hence Theorem \ref{5.6} and thus Propositions
 \ref{6.3}  are applicable. For the time being consider the case when $B$ is a germ of a smooth curve at a point $o$ and denote by $E$ the fiber in $X$ over $o$. 
 Since $\Phi |_E$ is a translation in Theorem \ref{0.2}  (3), for every $r \in R$ the preimage $\alpha^{-1}(r)$ is a curve. Thus, taking into consideration
 Corollary \ref{10.5}, we see that the assumptions of Proposition
 \ref{6.6} are valid. This implies that the polynomial curve $Z_{n-1}$ is smooth and we have $Z_{n-1} \simeq \C$  and $D_n\simeq \C^2$.
 Hence the polynomial curve $Z_{n-1}$ is smooth by Proposition \ref{6.6} and we have $Z_{n-1} \simeq \C$  and $D_n\simeq \C^2$. 
By assumption $D_n$ is also a smooth reduced fiber of the morphism $Q_n\to B$ and $Q_n \setminus D_n \simeq B^* \times \C^2$.
The Sathaye's theorem \cite{Sa} implies now that $Q_n$ is naturally isomorphic to $B \times \C^2$.  Similarly $X=Q_n \times \C$ when $B$ is a germ of a curve.

Let us return now to the general case when $B$ is a smooth affine curve. 
The local statement over germ of $B$ at $o$ implies the first global claim of Theorem \ref{0.2} (3) while the second one is the consequence of the fact (see \cite{Se})
that every affine manifold with a free $G_a$-action and an affine geometrical quotient is the direct product of that quotient and $\C$.
Thus Statement (3) of Theorem \ref{0.2} is valid while Statements (1) and (2) are true 
 by Corollaries \ref{5.7} and \ref{5.8}. Hence we are done with the proof of Theorem \ref{0.2}. \hfill    $ \square$
  
 \subsection{General case of an affine algebraic variety $X$ over a field $\bk$ of characteristic zero}\label{gc}  Recall that there is a one-to-one correspondence between the set of $G_a$-actions on an affine algebraic variety $X$ and 
 the set of locally nilpotent derivations (LNDs) on the algebra $\cB$ of regular functions on $X$ (e.g., see \cite{Fre}).
Suppose that $\nu$ is such an LND associated with a proper $G_a$-action $\Phi$ on $X$. Then condition that $\Phi$ is free (resp. a translation)
is equivalent to the fact that the ideal generated by $\nu (\cB )$ coincides with $\cB$ (resp. $\nu (\cB ) = \cB$). To show the validity of Theorem \ref{0.1}
for any $\bk$ of characteristic zero one needs to repeat the argument of Daigle from \cite[Theorem 3.2]{DaKa}.
 
Namely, if a fact is true for varieties over $\C$ it is true in the case when $X$ is  considered over a field which is a universal domain\footnote{Recall that a universal domain
is an algebraically closed field containing $\Q$ such that it has an infinite transcendence degree over $\Q$.} \cite{Ek}.
Thus consider a field extension $\bk'/\bk$ where $\bk'$ is a universal domain. Then $\nu$ extends to an LND $\nu'$ on $\cB' = \bk' \otimes_\bk B$
associated with a $G_a$-action $\Phi'$. Note that $\nu'$ is free since $\nu$ is. Similarly, $\Phi'$ is proper, since properness survives base extension \cite[Corollary 4.8]{Har}.
Thus under the assumptions of Theorem \ref{0.1} (with $\C$ replaced by $\bk'$) 
$\Phi'$ is a translation by the argument above, i.e. $\nu' : \cB' \to \cB'$ is surjective.  Since $\nu'$ is obtained by applying the functor $\bk' \otimes_\bk -$ to $\nu$
and since $\bk'$ is a faithfully flat $k$-module we see that $\nu : \cB \to \cB$ is also surjective. Therefore $\Phi$ is a translation and we have the following.

\bthm\label{7.0.1} Let $\bk$ be a field of characteristic zero and $\Phi$ be a proper $G_a$-action on the four-space $\AA_\bk^4$ preserving
a coordinate. Then $\Phi$ is a translation in a suitable polynomial coordinate system.

\ethm

%Note also that condition on triviality of second and third homology in Theorem \ref{0.2} is needed in the proof for one purpose only:  to guarantee that the restriction of $\Phi$
%to any fiber of $\varphi$ is a translation \cite{Ka04}. 
Similarly the argument before implies the following.

\bthm\label{7.0.2} Let the assumptions of Theorem \ref{0.2} hold with the only change that $X$ and $B$ are varieties over some field $\bk$ of characteristic zero
(and not over $\C$).  Suppose also that $\Phi$ is proper and the restriction of $\Phi$ to every fiber is a translation (as in Theorem \ref{0.2} (3)). Then $\Phi$ is a translation.

\ethm

\subsection{Conclusive Remark} 
In dimension 5 the analogue of Theorem \ref{0.1} is not true by the second Winkelmann's example in \cite{Wi}. This raises the question whether the properness
is the right condition for this type of problems. 
The author suspects that another condition may be more effective at least
in dimension 4 which leads to the following question.\\

{\em Suppose that the quotient of a free $G_a$-action $\Phi$ on $\C^4$ is isomorphic to $\C^3$. Is it true that $\Phi$ is a translation?}

\end{document}